\documentclass[a4paper, 11pt, twoside]{article}
\usepackage{mfpaperstuff}
\usepackage{needspace}
\usepackage{pgf}
\usepackage{tikz}
\usepackage[framemethod=tikz]{mdframed}

\newgeometry{margin=0.421in,includehead,includefoot}
\begin{document}
\newcommand\sq{1.5}
\newcommand\sssq{3.5}
\newcommand\ssssq{4.5}
\newcommand\hati{\hat\imath}
\newcommand\hatj{\hat\jmath}
\newcommand\roweq{\stackrel{\text{row}}\longleftrightarrow}
\newcommand\coleq{\stackrel{\text{column}}\longleftrightarrow}
\newcommand\ione{x_1}
\newcommand\ilast{x_s}
\newcommand\jone{y_1}
\newcommand\jtwo{y_2}
\newcommand\jlast{y_t}
\newcommand\ihat{\hat\imath}
\newcommand\po{{+}\negthinspace1}
\newcommand\nn[1]{N_{\doms#1}}
\newcommand\reg{^{\operatorname{reg}}}
\newcommand\swed[1]{\stackrel{#1}{\wed}}
\newcommand\inv[2]{\llbracket#1,#2\rrbracket}
\newcommand\one{\mathbbm{1}}
\newcommand\bsm{\begin{smallmatrix}}
\newcommand\esm{\end{smallmatrix}}
\newcommand{\rt}[1]{\rotatebox{90}{$#1$}}
\newcommand{\ol}{\overline}
\newcommand{\ul}{\underline}
\newcommand\partn{\mathcal{P}}
\newcommand{\sss}{\mathfrak{S}_}
\newcommand{\nin}{\notin}
\newcommand{\nchar}{\operatorname{char}}
\newcommand{\thmcite}[2]{\textup{\textbf{\cite[#2]{#1}}}\ }
\newcommand\zo{\bbn_0}
\newcommand{\sect}[1]{\section{#1}}
\newcommand\mmod{\ \operatorname{Mod}}
\newcommand\cgs\succcurlyeq
\newcommand\cls\preccurlyeq
\newcommand\UU{\mathcal{U}}
\newcommand\MM[1]{M^{\otimes#1}}
\newcommand\add{\operatorname{add}}
\newcommand\rem{\operatorname{rem}}
\newcommand\lra\longrightarrow
\newcommand\tru[1]{{#1}_-}
\newcommand\ste[1]{{#1}_+}
\newcommand\lad{\mathcal{L}}
\newcommand\hth[1]{\hat\Theta_{#1}}
\newcommand\ps[2]{\psi_{#1,#2}}
\renewcommand\hom{\operatorname{Hom}}
\newcommand\dpo{d\!\!+\!\!1}
\newcommand\bpt{b\!\!+\!\!2}
\newcommand\bph{b\!\!+\!\!3}
\newcommand\bpf{b\!\!+\!\!4}
\newcommand\bpv{b\!\!+\!\!5}
\newcommand\bpx{b\!\!+\!\!6}
\newcommand\vpo{v\!\!+\!\!1}
\newcommand\vpt{v\!\!+\!\!2}
\newcommand\xone{x_1}
\newcommand\xs{x_s}
\newcommand\zone{z_1}
\newcommand\ztwo{z_2}
\newcommand\zt{z_t}
\newcommand\ten{10}
\newcommand\eleven{11}
\newcommand\twelve{12}
\newcommand\ct\calu
\newcommand\semi[2]{\calto(#1,#2)}
\newcommand\ca[1]{\calt[#1]}
\newcommand\caltr{\calt_{\hspace{-2pt}\operatorname{r}}}
\newcommand\calto{\calt_{\hspace{-2pt}0}}
\newlength\frt
\newcommand\footy[1]{\text{\footnotesize$#1$}}
\newcommand\hdo[1]{\begin{tikzpicture}[baseline=0cm]\draw[thick,densely dotted](0,.13)--++(#1*.25,0);\end{tikzpicture}}
\newcommand\hds{\hdo1}
\newcommand\hdsq{\hdo{1.5}}
\newcommand\hdss{\hdo2}
\newcommand\hdsss{\hdo3}
\newcommand\hdsssq{\hdo{3.5}}
\newcommand\hdssss{\hdo4}
\newcommand\hdssssq{\hdo{4.5}}
\newcommand\hdsssss{\hdo5}
\newcommand\vdo[1]{\begin{tikzpicture}[baseline=0cm]\draw[thick,densely dotted](0,.13-#1*.125)--++(0,#1*.25);\end{tikzpicture}}
\newcommand\vds{\vdo1}
\newcommand\vdsq{\vdo{1.5}}
\newcommand\vdss{\vdo2}

\pdfpagewidth=12.5cm
\pdfpageheight=8.839cm
\setcounter{page}0
{\footnotesize This is the second author's version of a work that was accepted for publication in the Journal of Algebra. Changes resulting from the publishing process, such as peer review, editing, corrections, structural formatting, and other quality control mechanisms may not be reflected in this document. Changes may have been made to this work since it was submitted for publication. A definitive version was subsequently published in\\\textit{J.\ Algebra} \textup{\textbf{357} (2012) 235--62.\\http://dx.doi.org/10.1016/j.jalgebra.2012.01.035}\normalsize}
\newgeometry{a4paper,margin=1in,includehead,includefoot}
\title{Some new decomposable Specht modules}

\author{Craig J.\ Dodge\\
\normalsize Department of Mathematics, University at Buffalo, SUNY,\\
\normalsize 244 Mathematics Building, Buffalo, NY 14260, U.S.A.
\\\\
\large Matthew Fayers\\
\normalsize Queen Mary, University of London, Mile End Road, London E1 4NS, U.K.}

\toptitle

\markboth{Craig J.\ Dodge and Matthew Fayers}{Some new decomposable Specht modules}\pagestyle{myheadings}

\begin{abstract}
We present (with proof) a new family of decomposable Specht modules for the symmetric group in characteristic $2$. These Specht modules are labelled by partitions of the form $(a,3,1^b)$, and are the first new examples found for thirty years.  Our method of proof is to exhibit summands isomorphic to irreducible Specht modules, by constructing explicit homomorphisms between Specht modules.
\end{abstract}

\pdfpagewidth=8.27in
\pdfpageheight=11.69in

\renewcommand\baselinestretch{1.119}
\Yvcentermath1\Yboxdim{4.5mm}

\section{Introduction}

Let $n$ be a positive integer, and let $\sss n$ denote the symmetric group on $n$ letters. For any field $\bbf$, the \emph{Specht modules} form an important family of modules for $\bbf\sss n$. If $\bbf$ has characteristic zero, then the Specht modules are precisely the irreducible modules for $\bbf\sss n$. If $\bbf$ has positive characteristic, the simple $\bbf\sss n$-modules arise as quotients of certain Specht modules. In addition, the Specht modules arise as the `cell modules' for Murphy's cellular basis of $\bbf\sss n$.

A great deal of effort is devoted to determining the structure of Specht modules; in particular, finding the composition factors of Specht modules and the dimensions of the spaces of homomorphisms between Specht modules. In this paper, we consider the question of which Specht modules are decomposable. It is known that in odd characteristic the Specht modules are all indecomposable, so we can concentrate on the case where $\nchar(\bbf)=2$.  In fact, since any field is a splitting field for $\sss n$, we can assume that $\bbf=\bbf_2$. In this case, there are decomposable Specht modules, but remarkably few examples are known. Murphy \cite{gm} analysed the Specht modules labelled by `hook partitions', i.e.\ partitions of the form $(a,1^b)$, computing the endomorphism ring of every such Specht module (and thereby determining which ones are decomposable).  However, in the last thirty years no more progress seems to have been made.

Our main result is the discovery of a new family of decomposable Specht modules, the first examples of which were discovered by the two authors independently using computations with GAP and MAGMA. These new decomposable Specht modules are labelled by partitions of the form $(a,3,1^b)$, where $a,b$ are even. So in this paper we make a case study of partitions of this form; we are unable to apply Murphy's method to determine exactly which of these Specht modules are decomposable, but by considering homomorphisms between Specht modules, we are able to show which irreducible Specht modules arise as summands of these Specht modules. We then apply this result to determine which of our Specht modules have a summand isomorphic to an irreducible Specht module.

We now briefly indicate the layout of this paper.  In the next section, we recall some basic definitions and results in the representation theory of the symmetric group, which enable us to state our main results in Section \ref{resultsec}.  In Section \ref{homsec} we go into more detail on homomorphisms between Specht modules.  In Sections \ref{uvsec} and \ref{uv2sec} we consider the two classes of irreducible Specht modules which can occur as summands of our decomposable Specht modules.  We then apply these results in Section \ref{whichdec} to complete the proof of our main results.  Finally, we make some concluding remarks in Section \ref{concsec}.

\begin{acks}
The authors are indebted to David Hemmer, who first made us aware of each other's work and initiated this collaboration, and also invited the second author to SUNY Buffalo in September 2011, where some of this work was carried out.  This work continued during the `New York workshop on the symmetric group'; we are grateful to Rishi Nath of CUNY for inviting us to this conference.

The research of the first author was supported in part by NSA grant H98230-10-1-0192.
\end{acks}

\section{Background results}\label{backsec}

In this section, we summarise some basic results on the representation theory of the symmetric group. For brevity, we specialise some results to characteristic $2$, referring the reader to the literature for general results.

We begin by fixing a field $\bbf$; all our modules will be modules for the group algebra $\bbf\sss n$.  We assume familiarity with James's book \cite{j2}; in particular, we refer the reader there for the definitions of partitions, the dominance order, the permutation modules $M^\la$, the Specht modules $S^\la$ and the simple modules $D^\la$. We shall also briefly use the Nakayama Conjecture \cite[Theorem 21.11]{j2} which describes the block structure of the symmetric group.

We also need the following two results; recall that if $\la$ is a partition then $\la'$ denotes the conjugate partition.

\begin{lemma}\label{isospecht}
Suppose $\nchar(\bbf)=2$ and $\la$ is a partition such that $S^\la$ is irreducible. Then $S^\la\cong S^{\la'}$.
\end{lemma}

\begin{pf}
By \cite[Theorem 8.15]{j2} we have $S^\la\cong(S^{\la'})^\ast$, since the sign representation is trivial in characteristic $2$. But by \cite[Theorem 11.5]{j2}, all simple modules for the symmetric group are self-dual.
\end{pf}

\begin{lemma}\label{815hom}
If $\la,\mu$ are partitions of $n$, then
\[
\dim_\bbf\hom_{\bbf\sss n}(S^\la,S^\mu)=\dim_\bbf\hom_{\bbf\sss n}(S^{\mu'},S^{\la'}).
\]
\end{lemma}

\begin{pf}
This also follows from \cite[Theorem 8.15]{j2}.
\end{pf}

\subsection{Regularisation}

We recall here a useful lemma which we shall use later; this is due to James, although it does not appear in the book \cite{j2}.  We concentrate on the special case where $\bbf$ has characteristic~$2$, referring to \cite{j1} for the full result.

For any $l\gs1$, the $l$th \emph{ladder} in $\bbn^2$ is
\[
\call_l=\lset{(i,j)}{i+j=l+1}.
\]
If $\la$ is a partition, the \emph{$2$-regularisation} of $\la$ is the partition $\la\reg$ whose Young diagram is obtained by moving the nodes in $[\la]$ as high as possible within their ladders.  For example, $(8,3,1^6)\reg=(8,7,2)$, as we see from the following Young diagrams, in which nodes are labelled according to the ladders in which they lie.
\Yvcentermath0
\[
\young(12345678,234,3,4,5,6,7,8)\qquad\qquad
\young(12345678,2345678,34)
\]
\Yvcentermath1
It is a simple exercise to show that $\la\reg$ is a $2$-regular partition, and we have the following result.

\begin{thm}\thmcite{j1}{Theorem A}\label{jreg}
Suppose $\la$ and $\mu$ are partitions of $n$, with $\mu$ $2$-regular.  Then $[S^\la:D^{\la\reg}]=1$, while $[S^\la:D^\mu]=0$ if $\mu\ndom\la\reg$.
\end{thm}

In this paper we shall be concerned with the Specht modules labelled by partitions of the form $(a,3,1^b)$; so we compute the regularisations of these partitions.

\begin{lemma}\label{reg}
Suppose $a\gs4$ and $b\gs2$. Then
\[
(a,3,1^b)\reg=
\begin{cases}
(a,b+1,2)&(a>b)\\
(b+2,a-1,2)&(a\ls b).
\end{cases}
\]
\end{lemma}

\subsection{Irreducible Specht modules}

It will be very helpful to know the classification of irreducible Specht modules, which (in characteristic $2$) was discovered by James and Mathas \cite{jmp2}.  If $k$ is a non-negative integer we let $l(k)$ denote the smallest positive integer such that $2^{l(k)}>k$.

\begin{thm}\thmcite{jmp2}{Main Theorem}\label{irrspecht}
Suppose $\mu$ is a partition of $n$ and $\nchar(\bbf)=2$.  Then $S^\mu$ is irreducible if and only if one of the following occurs:
\begin{itemize}
\item
$\mu_i-\mu_{i+1}\equiv-1\ppmod{2^{l(\mu_{i+1}-\mu_{i+2})}}$ for each $i\gs1$;
\item
$\mu'_i-\mu'_{i+1}\equiv-1\ppmod{2^{l(\mu'_{i+1}-\mu'_{i+2})}}$ for each $i\gs1$;
\item
$\mu=(2^2)$.
\end{itemize}
\end{thm}

Note that $\mu$ satisfies the first condition in the theorem if and only if $\mu'$ satisfies the second.  In view of Lemma \ref{isospecht} (and since we shall only be considering values of $n$ greater than $4$) we may assume that any irreducible Specht module is of the form $S^\mu$ where $\mu$ satisfies the first condition in the theorem.

\section{The main results}\label{resultsec}

In this section, we describe the new family of decomposable Specht modules discussed in this paper, and the method we use to prove decomposability. \emph{For the rest of this paper, we assume that $\bbf$ has characteristic $2$.}

Computer calculations show that the first few decomposable Specht modules which are not labelled by hook partitions have labelling partitions of the form $(a,3,1^b)$ (and their conjugates) with $a,b$ even positive integers. So in this paper we make a case study of this family of partitions. Our technique is different from that of Murphy \cite{gm}, and is weaker in the sense that we cannot always when tell for certain whether one of our Specht modules is decomposable. However, in the cases where we can show decomposability, we have the advantage of being able to describe one summand explicitly.

More specifically, our main result is a determination of exactly which irreducible Specht modules occur as summands of the Specht modules $S^{(a,3,1^b)}$. The technique we use is to consider homomorphisms between Specht modules, and the set-up for computing such homomorphisms is described in Section \ref{homsec}.

We use homomorphisms between Specht modules in the following way. Suppose $\la,\mu$ are partitions of $n$.  Then it is a straightforward result that $S^\mu$ occurs as a summand of $S^\la$ if and only if there are homomorphisms $\gamma:S^\mu\to S^\la$ and $\delta:S^\la\to S^\mu$ such that $\delta\circ\gamma$ is the identity on $S^\mu$.  If we assume in addition that $S^\mu$ is irreducible, then by Schur's Lemma we just need to show that $\delta\circ\gamma$ is non-zero.

Some effort has been devoted to computing the space of homomorphisms between two Specht modules, beginning with the paper of the second author and Martin \cite{fm}. In fact, there is now an explicit algorithm which computes the homomorphism space $\hom_{\bbf\sss n}(S^\la,S^\mu)$ except when $\nchar(\bbf)=2$ and $\la$ is $2$-singular. Even in this exceptional case, this technique can be used to construct some homomorphisms between $S^\la$ and $S^\mu$, though only if $\la$ dominates $\mu$.

In our situation, the partitions $\mu$ we shall consider are always $2$-regular, because (as long as $n\neq4$) every irreducible Specht module in characteristic $2$ has the form $S^\mu$ for $\mu$ $2$-regular.  So we are able to compute the space $\hom_{\bbf\sss n}(S^\mu,S^\la)$.  Computing the space $\hom_{\bbf\sss n}(S^\la,S^\mu)$ is harder, because $\la$ is $2$-singular, so we fall foul of the exception above.  To circumvent this, we use Lemma \ref{isospecht}, which allows us to take $\delta$ to be an element of $\hom_{\bbf\sss n}(S^\la,S^{\mu'})$. By Lemma \ref{815hom} this has the same dimension as $\hom_{\bbf\sss n}(S^\mu,S^{\la'})$, which we can compute because $\mu$ is $2$-regular.  Having established the dimension of this space, we can construct all possible homomorphisms $\delta$, and then check the condition $\delta\circ\gamma\neq0$.

In this way, we can find all summands of $S^\la$ which are isomorphic to irreducible Specht modules $S^\mu$. In fact, we can restrict attention to a small set of candidate Specht modules $S^\mu$, as follows. Assuming $S^\mu$ is irreducible and $\mu$ is $2$-regular, $S^\mu$ is isomorphic to the simple module $D^\mu$; therefore in order for $S^\mu$ to appear as a summand of $S^\la$, the decomposition number $[S^\la:D^\mu]$ must be non-zero.  Therefore by Theorem \ref{jreg}, $\mu$ must dominate $\la\reg$, and by Lemma \ref{reg} this has the form $(x,y,2)$ for some $x,y$.  So we may assume that $\mu$ has the form $(u,v,w)$ for some $u>v>w\ls2$. Furthermore, $S^\mu$ and $S^\la$ must lie in the same block of $\bbf\sss n$. Using the Nakayama Conjecture, this means that $u,w$ must be even, while $v$ is odd. So we can restrict attention to $\mu$ of the form $(u,v)$ or $(u,v,2)$ where $u$ is even (and hence $v$ is odd).  We apply the technique described above to these two types of partition in Sections \ref{uvsec} and 
\ref{uv2sec}.  Our main result is the following.

\begin{mdframed}[innerleftmargin=3pt,innerrightmargin=3pt,innertopmargin=0pt,innerbottommargin=3pt,roundcorner=5pt,innermargin=-3pt,outermargin=-3pt]
\begin{thm}\label{main}
Suppose $\la=(a,3,1^b)$ is a partition of $n$, where $a,b$ are positive even integers with $a\gs4$, and suppose $\mu$ is a partition of $n$ such that $S^\mu$ is irreducible.  Then $S^\la$ has a direct summand isomorphic to $S^\mu$ if and only if one of the following occurs.
\begin{enumerate}
\item
$\mu$ or $\mu'$ equals $(u,v)$, where $v\equiv3\ppmod4$ and $\mbinom{u-v}{a-v}$ is odd.
\item
$\mu$ or $\mu'$ equals $(u,v,2)$, where $\mbinom{u-v}{a-v}$ is odd.
\end{enumerate}
\end{thm}
\end{mdframed}

\vspace{\topsep}
Using this result, we can show that most of the Specht modules under consideration are decomposable. Specifically, we have the following result.
\vspace{\topsep}

\begin{mdframed}[innerleftmargin=3pt,innerrightmargin=3pt,innertopmargin=0pt,innerbottommargin=3pt,roundcorner=5pt,innermargin=-3pt,outermargin=-3pt]
\begin{cory}\label{maincor}
Suppose $a,b$ are positive even integers with $a\gs4$, and let $\la=(a,3,1^b)$. Then $S^\la$ has a summand isomorphic to an irreducible Specht module if and only if at least one of the following occurs:
\begin{itemize}
\item
$a+b\equiv0$ or $2\pmod 8$, $a\gs6$ and $b\gs4$;
\item
$a+b\equiv2\ppmod4$ and $\mbinom{a+b-3}{a-3}$ is odd;
\item
$a+b\equiv0\ppmod4$ and $\mbinom{a+b-9}{a-5}$ is odd.
\end{itemize}
\end{cory}
\end{mdframed}

\section{Computing the space of homomorphisms between two Specht modules}\label{homsec}

In this section, we explain the set-up for computing the space of homomorphisms between two Specht modules.  We begin with a revision of some material from \cite{j2}, before citing some results of the second author and Martin.

\subsection{Homomorphisms from Specht modules to permutation modules}

Suppose $\mu$ and $\la$ are partitions of $n$. Since $S^\la\ls M^\la$, any homomorphism from $S^\mu$ to $S^\la$ can be regarded as a homomorphism from $S^\mu$ to $M^\la$.  This is very useful, because if $\mu$ is $2$-regular (or if $\nchar(\bbf)\neq2$), then the space $\hom_{\bbf\sss n}(S^\mu,M^\la)$ can be described explicitly. Furthermore, using the Kernel Intersection Theorem below, one can check whether the image of a homomorphism $\theta:S^\mu\to M^\la$ lies in $S^\la$.

We now make some more precise definitions.  We take $\la,\mu$ as above, but we now allow $\la$ to be any composition of $n$, not necessarily a partition. A \emph{$\mu$-tableau of type $\la$} is a function $T$ from the Young diagram $[\mu]$ to $\bbn$ with the property that for each $i\in\bbn$ there are exactly $\la_i$ nodes of $[\mu]$ mapped to $i$. Such a tableau is usually represented by drawing $[\la]$ with a box for each node $\fkn$, filled with the integer $T(\fkn)$. $T$ is \emph{row-standard} if the entries in this diagram are weakly increasing along the rows, and is \emph{semistandard} if the entries are weakly increasing along the rows and strictly increasing down the columns.

We write $\caltr(\mu,\la)$ for the set of row-standard $\mu$-tableaux of type $\la$, and $\calto(\mu,\la)$ for the set of semistandard $\mu$-tableaux of type $\la$.  For each $T\in\caltr(\mu,\la)$, James defines a homomorphism $\Theta_T:M^\mu\to M^\la$ (over any field), whose precise definition we do not need here.  The restriction of $\Theta_T$ to $S^\mu$ is denoted $\hth T$.  Now we have the following.

\begin{thm}\thmcite{j2}{Lemma 13.11 \& Theorem 13.13}\label{semi}
The set
\[
\lset{\hth T}{T\in\calto(\mu,\la)}
\]
is linearly independent, and spans $\hom_{\bbf\sss n}(S^\mu,M^\la)$ if $\mu$ is $2$-regular.
\end{thm}

\subsection{The Kernel Intersection Theorem}

Now return to the assumption that $\la$ is a partition. As a consequence of Theorem \ref{semi}, in order to compute $\hom_{\bbf\sss n}(S^\mu,S^\la)$ when $\mu$ is $2$-regular, we just need to find all linear combinations $\theta$ of the homomorphisms $\hth T$ for $T\in\calto(\mu,\la)$ for which the image of $\theta$ lies in $S^\la$. Even when $\mu$ is not $2$-regular, homomorphisms from $S^\mu$ to $S^\la$ can very often be expressed in this way.

In order to determine whether the image of such a homomorphism $\theta$ lies in $S^\la$, we use another theorem of James which provides an alternative definition of $S^\la$. For any pair $(d,t)$ with $d\gs 1$ and $1\ls t\ls \la_{d+1}$, there is a homomorphism $\ps dt:M^\la\to M^\nu$, where $\nu$ is a composition depending on $\la,d,t$. Again, we refer the reader to \cite[\S17]{j2} for the definition; we warn the reader that the homomorphism $\ps dt$ is called $\psi_{d,\la_{d+1}-t}$ in [\textit{loc.\ cit.}]. The importance of the homomorphisms $\ps dt$ is in the following.
\begin{thm}\thmlabel{The Kernel Intersection Theorem}\thmcite{j2}{Corollary 17.18}\label{kit}
Suppose $\la$ is a partition of $n$. Then
\[
S^\la=\bigcap_{\begin{smallmatrix}d\gs1\\1\ls t\ls\la_{d+1}\end{smallmatrix}}\ker(\ps dt).
\]
\end{thm}

This provides a clear strategy for computing $\hom_{\bbf\sss n}(S^\mu,S^\la)$: find all linear combinations $\theta$ of the homomorphisms $\hth T$ such that $\ps dt\circ\theta=0$ for every $d,t$. Fortunately, it is known how to compute the composition $\ps dt\circ\hth T$ when $T\in \caltr(\mu,\la)$. For our next few results, we need to introduce some more notation. For any multiset $A$ of positive integers, let $A_i$ denote the number of $i$s in $A$. If $A,B$ are multisets, we write $A\sqcup B$ for the multiset with $(A\sqcup B)_i=A_i+B_i$ for all $i$. Given a row-standard tableau $T$, we write $T^j$ for the multiset of entries in row $j$ of $T$.

\begin{lemma}\thmcite{fm}{Lemma 5}\label{lemma5}
Suppose $\la,\mu$ are partitions of $n$, $T\in\caltr(\mu,\la)$, $d\in\bbn$ and $1\ls t\ls\la_{d+1}$. Let $\cals$ be the set of all row-standard tableaux which can be obtained from $T$ by replacing $t$ of the entries equal to $d+1$ in $T$ with $d$s.  Then
\[
\ps dt\circ\Theta_T=\sum_{S\in\cals}\prod_{j\gs1}\binom{S^j_d}{T^j_d}\Theta_S.
\]
\end{lemma}

The slight difficulty with using this lemma to compute homomorphism spaces is that the tableaux $S$ in the lemma are not always semistandard; so it can be difficult to tell whether a particular linear combination is zero when restricted to $S^\mu$. To circumvent this, we recall another lemma from \cite{fm} which gives certain linear relations between the homomorphisms $\hth T$, and often enables us to write a homomorphism $\hth T$ in terms of semistandard homomorphisms.

\begin{lemma}\thmcite{fm}{Lemma 7}\label{lemma7}
Suppose $\mu$ is a partition of $n$ and $\la$ a composition of $n$, and $i,j,k$ are positive integers with $j\neq k$ and $\mu_j\gs\mu_k$. Suppose $T\in\caltr(\mu,\la)$, and let $\cals$ be the set of all $S\in\caltr(\mu,\la)$ such that:
\begin{itemize}
\item
$S^j_i=T^j_i+T^k_i$;
\item
$S^j_l\ls T^j_l$ for every $l\neq i$;
\item
$S^l=T^l$ for all $l\neq j,k$.
\end{itemize}
Then
\[\hth T = (-1)^{T^k_i}\sum_{S\in\cals}\prod_{l\gs1}\binom{S^k_l}{T^k_l}\hth S.\]
\end{lemma}

Informally, a tableau in $\cals$ is a tableau obtained from $T$ by moving all the $i$s from row $k$ to row $j$, and moving some multiset of entries different from $i$ from row $j$ to row $k$.

One very simple case of Lemma \ref{lemma7} which we shall apply frequently is the following: if $T\in\caltr(\mu,\la)$ and for some $i,j,k$ we have $T^j_i+T^k_i>\max\{\mu_j,\mu_k\}$, then $\hth T=0$.

Lemma \ref{lemma7} turns out to be very useful for expressing a tableau homomorphism in terms of semistandard homomorphisms.  However, we shall occasionally need to use the following alternative.

\begin{lemma}\label{newsemilem}
Suppose $\la$ and $\mu$ are partitions of $n$, and $T$ is a row-standard $\la$-tableau of type $\mu$.  Suppose $i\gs1$, and $A,B,C$ are multisets of positive integers such that $|B|>\la_i$ and $A\sqcup B\sqcup C=T^i+T^{i+1}$.  Let $\calb$ be the set of all pairs $(D,E)$ of multisets such that $|D|=\la_i-|A|$ and $B=D\sqcup E$.  For each such pair $(D,E)$, define $T_{D,E}$ to be the row-standard tableau with
\[
T_{D,E}^j=
\begin{cases}
A\sqcup D&(j=i)\\
C\sqcup E&(j=i+1)\\
T^j&(\text{otherwise}).
\end{cases}
\]
Then
\[
\sum_{(D,E)\in\calb}\prod_{i\gs1}\binom{A_i+D_i}{D_i}\binom{C_i+E_i}{E_i}\hth{T_{D,E}}=0.
\]
\end{lemma}

This lemma appears in the second author's forthcoming paper \cite{garnir} where it is proved in the wider context of Iwahori--Hecke algebras; however, a considerably easier proof exists in the symmetric group case.  In \cite{garnir}, Lemma \ref{newsemilem} is used to provide an explicit fast algorithm for writing a tableau homomorphism as a linear combination of semistandard homomorphisms.

We now give another result which will help us in showing that a linear combination of row-standard homomorphisms is non-zero without having to go through the full process of expressing it as a linear combination of semistandard homomorphisms.  This concerns the \emph{dominance order} on tableaux.

Suppose $\mu$ is a partition, and $S,T$ are row-standard $\mu$-tableaux of the same type.  We say that $S$ dominates $T$, and write $S\dom T$, if it is possible to  get from $S$ to $T$ by repeatedly swapping an entry of $S$ with a larger entry in a lower row (and re-ordering with each row).  We warn the reader that this is not quite the same as the dominance order described in \cite[13.8]{j2}.  For example, the dominance order on $\caltr\left((3,2),(2^2,1)\right)$ is given by the following Hasse diagram.

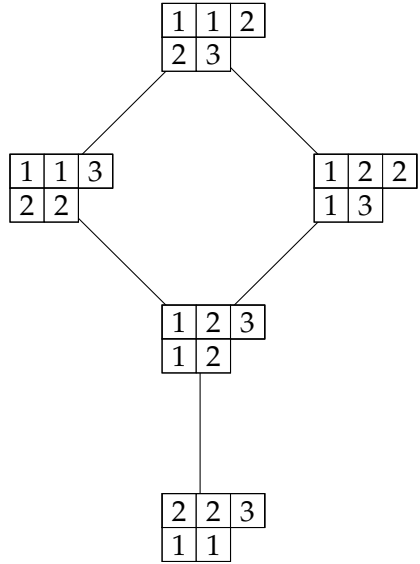
\[
\begin{tikzpicture}
\draw(0.5,0.5)--++(0,2.5)--++(2,2)--++(-2,2)--++(-2,-2)--++(2,-2);
\tyoung(0cm,0.5cm,223,11)
\tyoung(0cm,3cm,123,12)
\tyoung(2cm,5cm,122,13)
\tyoung(-2cm,5cm,113,22)
\tyoung(0cm,7cm,112,23)
\end{tikzpicture}
\]

Now we have the following lemma.

\begin{lemma}\label{semidom}
Suppose $\mu$ is a partition of $n$ and $\la$ a composition of $n$, and $T\in\caltr(\mu,\la)$.  If we write
\[
\hth T=\sum_{S\in\calto(\mu,\la)}a_S\hth S,
\]
then $a_S\neq0$ only if $S\dom T$.
\end{lemma}

\begin{pf}
For this proof, we adopt the set-up of \cite[\S13]{j2}. If we fix a $\mu$-tableau $t$ (of type $(1^n)$), then $t$ determines a natural bijection between the set of $\la$-tabloids and the set of $\mu$-tableaux of type $\la$.  We identify tabloids with tableaux according to this bijection.  We also let $e_t$ denote the polytabloid indexed by $t$, which generates $S^\mu$. Given $\la$-tableaux $U,V$, we write $U\roweq V$ if $V$ can be obtained from $U$ by permuting entries within rows, and we define $U\coleq V$ similarly.  Now given $U\in\caltr(\mu,\la)$, we have
\[
\hth U(e_t)=\sum_{(V,W)}W,
\]
where we sum over all pairs $(V,W)$ of $\la$-tableaux such that $U\roweq V\coleq W$ and $V$ has distinct entries in each column. (In general there are also signs determined by the column permutations, but in characteristic $2$ we can neglect these.)

Now if $S$ is a semistandard tableau with $U\roweq V\coleq S$ for some $V$, then (since the entries in each column of $S$ are increasing) we must have $S\dom U$.  On the other hand, if $U$ is semistandard, then the only tableau $V$ such that $U\roweq V\coleq U$ is $U$ itself.  Hence for any semistandard $U$, the coefficient of $U$ in $\hth U(e_t)$ is $1$.

Now suppose $S_0\in\calto(\mu,\la)$ is  such that $a_{S_0}\neq0$ and $S_0$ is minimal (with respect to $\dom$) with this property.  Then the coefficient of $S_0$ in
\[
\hth T(e_t)=\sum_{S\in\calto(\mu,\la)}a_S\hth S(e_t)
\]
is $a_{S_0}$.  So the coefficient of $S_0$ in $\hth T(e_t)$ is non-zero, and hence $S_0\dom T$.  Any $S\in\calto(\mu,\la)$ for which $a_S\neq0$ dominates some such minimal tableau $S_0$, and so dominates $T$.
\end{pf}

Using the results in this section, it will be possible to compute $\hom_{\bbf\sss n}(S^\mu,S^\la)$ in the cases of interest to us. We remark that it is often easier to express such homomorphisms as linear combinations of non-semistandard homomorphisms; in particular, the conditions $\ps dt\circ\theta=0$ can be easier to check.  Of course, when doing this we have to be careful to show that the homomorphisms we construct are non-zero.

\subsection{Composition of tableau homomorphisms}

It will also be important in this paper to compute compositions of homomorphisms between Specht modules.  It is well understood how to compose two tableau homomorphisms; indeed, computing this composition is the same as computing the structure constants for the Schur algebra.  We give this result, of which Lemma \ref{lemma5} is a special case.  This result is easy to prove and well known (indeed, `quantised' versions appear in the literature) but it does not seem to appear explicitly.  However, translating to the language of the Schur algebra (where $\Theta_T$ corresponds to a basis element $\xi_{i,j}$) it amounts to the Multiplication rule (2.3b) in Green's monograph \cite{gr}.

Recall that if $S$ is a tableau, then $S^j$ denotes the multiset of entries in row $j$ of $S$, and in particular $S^j_i$ denotes the number of entries equal to $i$ in row $j$ of $S$. If $x_1,x_2,\dots$ are non-negative integers with finite sum $x$, we write $\mbinom{x}{x_1,x_2,\dots}$ for the corresponding multinomial coefficient.

\begin{propn}\label{tabcomp}
Suppose $\la,\mu,\nu$ are compositions of $n$, $S$ is a $\la$-tableau of type $\mu$ and $T$ is a $\mu$-tableau of type $\nu$. Let $\calx$ be the set of all collections $X=(X^{ij})_{i,j\gs1}$ of multisets such that
\[
|X^{ij}|=S^j_i\quad\text{ for each $i$, $j$,}\qquad\bigsqcup_{j\gs1}X^{ij}=T^i\quad\text{ for each $i$.}
\]
For $X\in\calx$, let $U_X$ denote the row-standard $\la$-tableau with $(U_X)^j=\bigsqcup_{i\gs1}X^{ij}$.  Then
\[
\Theta_T\circ\Theta_S = \sum_{X\in\calx}\prod_{i,j\gs1}\binom{X^{1j}_i+X^{2j}_i+X^{3j}_i+\dots}{X^{1j}_i,X^{2j}_i,X^{3j}_i,\dots}\Theta_{U_X}.
\]
\end{propn}

\section{Irreducible summands of the form $S^{(u,v)}$}\label{uvsec}

In this section, we find all cases where one of our Specht modules $S^{(a,3,1^b)}$ has a summand isomorphic to an irreducible Specht module of the form $S^{(u,v)}$, where $u$ is even and $v$ is odd.  Throughout, we continue to assume that $a,b$ are positive even integers with $a\gs4$, and we let $n=a+b+3$.  By Theorem \ref{jreg} and Lemma \ref{reg}, $D^{(u,v)}$ cannot appear as a composition factor of $S^{(a,3,1^b)}$ unless $(u,v)\dom(a,3,1^b)\reg$, which is the partition $(\max\{a,b+2\},\min\{a-1,b+1\},2)$.  So we may assume that this is the case, which is the same as saying $v\ls\min\{a+1,b+3\}$.  For easy reference, we set out notation and assumptions for this section.

\needspace{10em}
\begin{mdframed}[innerleftmargin=3pt,innerrightmargin=3pt,innertopmargin=3pt,innerbottommargin=3pt,roundcorner=5pt,innermargin=-3pt,outermargin=-3pt]
\noindent\textbf{Assumptions and notation in force throughout Section \ref{uvsec}:}

$\la=(a,3,1^b)$ and $\mu=(u,v)$, where $a,b,u,v$ are positive integers with $a,b,u$ even, $a\gs4$, $u>v$, $n=a+b+3=u+v$ and $v\ls\min\{a+1,b+3\}$.
\end{mdframed}

\subsection{Homomorphisms from $S^\la$ to $S^{\mu'}$}\label{hlamu'1}

In this subsection we consider $\bbf\sss n$-homomorphisms from $S^\la$ to $S^{\mu'}$. We begin by constructing such a homomorphism in the case where $3\ls v\ls a-1$.

Let $\ct$ be the set of $\la$-tableaux having the form
\[
\gyoung(;1;2;3_\sq\hdsq;v;\star_\sq\hdsq;\star,;1;\star;\star,;\star,|2\vdss,;\star)
\]
in which the $\star$s represent the numbers from $2$ to $u$, and in which
\begin{itemize}
\item
the entries along each row are strictly increasing,
\item
the entries down each column are weakly increasing.
\end{itemize}
Now define
\[
\sigma=\sum_{T\in\ct}\hth T.
\]

\begin{propn}\label{sigmahom}
With the assumptions and notation above, we have $\ps dt\circ\sigma=0$ for each $d,t$.
\end{propn}

\begin{pf}
First take $v<d\ls u$ and $t=1$.  If $T\in\ct$, then $T$ contains a single $d$ and a single $d+1$. If these occur in the same row or the same column of $T$, then $\ps d1\circ\hth T=0$ by Lemma \ref{lemma5} and Lemma \ref{lemma7}. Otherwise, there is another tableau $T'\in\ct$ obtained by interchanging the $d$ and the $d+1$. By Lemma \ref{lemma5} we have $\ps d1\circ(\hth T+\hth{T'})=0$.  Hence by summing $\ps d1\circ\hth T$ over all $T\in\ct$, we get zero.  A similar argument applies in the case $d=v$.

If $1\ls d<v$ and $t=2$, then we have $\ps dt\circ\hth T=0$ for each $T\in\ct$ just using Lemma \ref{lemma5}. Now take $2\ls d<v$ and $t=1$, and consider a tableau $T\in\ct$. There are a single $d$ and a single $d+1$ below row $1$. If these lie in the same row or column, then $\ps d1\circ\hth T=0$.  Otherwise, let $T'$ be the tableau obtained by interchanging the $d$ and the $d+1$ below row $1$. Then $\ps d1\circ(\hth T+\hth{T'})$, and we are done.

We are left with the case $d=t=1$. Applying Lemma \ref{lemma5}, we find that $\ps11\circ\theta$ is the sum of homomorphisms labelled by tableaux
\[
\gyoung(;1;2;3_\sq\hdsq;v;\star_\sq\hdsq;\star,;1;\star;\star,;1,;\star,|2\vdss,;\star)
\]
in which the $\star$s now represent the numbers from $3$ to $u$, and where the entries are strictly increasing along rows and weakly increasing down columns. Now we apply Lemma \ref{lemma7} to each of these homomorphisms to move the $1$ from row $3$ to row $2$, and then to reorder rows $3,\dots,b+2$. We obtain a sum of tableaux of the form
\[
\gyoung(;1;2;3_\sq\hdsq;v;\star_\sq\hdsq;\star,;1;1;\star,;\star,|2\vdss,;\star),
\]
but each tableau occurs $b$ times in this way. Since $b$ is even, we have zero.
\end{pf}
Now we need to check that $\sigma\neq0$, which is not obvious because the tableaux involved are not semistandard.

\begin{propn}\label{sigmanz}
With the notation above, $\sigma\neq0$.
\end{propn}

\begin{rmk}
The version of this paper published in the Journal of Algebra includes a fallacious proof of Proposition \ref{sigmanz}; the proof below replaces it. The authors are grateful to Sin\'ead Lyle for pointing out the error.
\end{rmk}

\begin{pf}
We'll use Lemma \ref{semidom}. Consider the semistandard tableau
\[
S=\footy{\gyoungx(1.2,;1;1;2_2\hdss;v;\bpv;\bpx_2\hdss;u,;2;\bph;\bpf,;3,;4,|2\vdss,;\bpt)}.
\]
We'll show that when $\sigma$ is expressed as a linear combination of semistandard homomorphisms, $\hth S$ occurs with non-zero coefficient, and hence $\sigma\neq0$.  Given $T\in\ct$, consider expressing $\hth T$ as a linear combination of semistandard homomorphisms.  By Lemma \ref{semidom}, $\hth S$ can only occur if $S\dom T$; so we can ignore all $T\in\ct$ for which $S\ndom T$. In particular, we need only consider those tableaux in $\ct$ which have $b+5,\dots,u$ in the first row and $b+3,b+4$ in the top two rows. If we assume for the moment that $v<b+3$, then the tableaux $T\in\ct$ that we need to consider are those of the following forms:
{\allowdisplaybreaks\begin{align*}
T[i]&=\footy{\gyoungx(1.2,;1;2;3_2\hdss;v;i;\bpv;\bpx_2\hdss;u,;1;\bph;\bpf,;2,;\vds,;v,;\vpo,;\vds,;\hati,;\vds,;\bpt)}\tag*{for $v<i\ls b+2$;}\\
U[i]&=\footy{\gyoungx(1.2,;1;2;3_2\hdss;v;\bph;\bpv;\bpx_2\hdss;u,;1;i;\bpf,;2,|\sq\vdsq,;\hati,|\sq\vdsq,;\bpt)}\tag*{for $2\ls i\ls b+2$;}\\
V[i]&=\footy{\gyoungx(1.2,;1;2;3_2\hdss;v;\bpf;\bpv;\bpx_2\hdss;u,;1;i;\bph,;2,|\sq\vdsq,;\hati,|\sq\vdsq,;\bpt)}\tag*{for $2\ls i\ls b+2$.}
\end{align*}}
As usual, the $\hati$ in the first column indicates that $i$ does not appear in that column.

First consider the tableau $T[i]$, and apply Lemma \ref{lemma7} to move the $1$ from row $2$ to row $1$. Of the tableaux appearing in the resulting expression, the only ones dominated by $S$ are
{\allowdisplaybreaks\begin{align*}
T'[i]&=\footy{\gyoungx(1.2,;1;1;2;3_2\hdss;v;\bpv;\bpx_2\hdss;u,;i;\bph;\bpf,;2,;\vds,;v,;\vpo,;\vds,;\hati,;\vds,;\bpt)}\\
\intertext{and the tableaux}
T'[i,j]&=\footy{\gyoungx(1.2,;1;1;2;\hds;\hatj;\hds;v;i;\bpv;\bpx_2\hdss;u,;j;\bph;\bpf,;2,;\vds,;v,;\vpo,;\vds,;\hati,;\vds,;\bpt)}\tag*{for $2\ls j\ls v$.}\\
\intertext{Applying Lemma \ref{lemma7} to $T'[i]$ to move the $2$ from row $3$ to row $2$, we obtain three tableaux, but two of these are not dominated by $S$. The other one is}
T''[i]&=\footy{\gyoungx(1.2,;1;1;2;3_2\hdss;v;\bpv;\bpx_2\hdss;u,;2;\bph;\bpf,;i,;3,;\vds,;v,;\vpo,;\vds,;\hati,;\vds,;\bpt)},
\end{align*}}
and $i-3$ more applications of Lemma \ref{lemma7} show that $\hth{T''[i]}=\hth S$.

Now consider applying Lemma \ref{lemma7} to $T'[i,j]$, to move the $2$ from row $3$ to row $2$. If $j=2$, then neither of the two tableaux obtained is dominated by $S$.  If $j>2$, then two of the three tableaux obtained are not dominated by $S$; the third has rows $3$ and $j+1$ both equal to $\young(j)$, so the resulting homomorphism is zero by Lemma \ref{lemma7}.

So we conclude that $\hth{T[i]}$ equals $\hth S$ plus a linear combination of homomorphisms indexed by tableaux not dominated by $S$.\\\hspace*{17pt}Now we consider applying Lemma \ref{lemma7} to $U[i]$, moving the $1$ up from row $2$. The tableaux obtained that are dominated by $S$ are $T'[i]$ and the tableaux
{\allowdisplaybreaks\begin{align*}
U'[i,j]&=\footy{\gyoungx(1.2,;1;1;2;\hds;\hatj;\hds;v;\bph;\bpv;\bpx_2\hdss;u,;j;i;\bpf,;2,;\vds,;\hati,;\vds,;\bpt)}\tag*{for $2\ls j\ls v$ with $i\neq j$.}\\
\intertext{\parbox{\linewidth}{(Note that if $i<j$ then $i$ and $j$ in the second row should be written the other way round; the case $i=j$ does not occur because the accompanying coefficient would be $\binom21=0$.)\\\hspace*{17pt}If $i=2$, then $U'[i,j]$ is a semistandard tableau different from $S$. If $i>2$ then we apply Lemma \ref{lemma7} to move the $2$ from row $3$ to row $2$; neglecting the tableau not dominated by $S$ and (in the case $j>2$) neglecting the tableau with two rows equal to $\young(j)$, the only tableau we get is}}
U''[i,j]&=\footy{\gyoungx(1.2,;1;1;2;\hds;\hatj;\hds;v;\bph;\bpv;\bpx_2\hdss;u,;2;j;\bpf,;i,;3,;\vds,;\hati,;\vds,;\bpt)};
\end{align*}}
$i-3$ more applications of Lemma \ref{lemma7} show that $\hth{U''[i,j]}$ equals a semistandard homomorphism different from $\hth S$.

We conclude that $\hth{U[i]}$ equals $\hth S$ plus a linear combination of homomorphisms indexed by tableaux which are either not dominated by $S$ or semistandard and different from $S$. The homomorphism $\hth{V[i]}$ is analysed in exactly the same way, interchanging $b+3$ and $b+4$.

Putting these cases together, we find that the coefficient of $\hth S$ in $\sigma$ is the total number of tableaux of the form $T[i]$, $U[i]$ or $V[i]$, i.e.\ $(b+2-v)+2(b+1)$, which is odd.

It remains to consider the case $v=b+3$. In this case only the tableaux $V[i]$ appear, but the analysis of these tableaux is exactly the same, so the coefficient of $\hth S$ in $\sigma$ is the number of tableaux $V[i]$, i.e.\ $b+1$, which again is odd.
\end{pf}

It turns out that up to scaling, $\sigma$ is the only homomorphism from $S^\la$ to $S^\mu$.

\begin{propn}\label{cdhomdim1}
With $\la,\mu$ as above,
\[
\dim_\bbf\hom_{\bbf\sss n}(S^\la,S^{\mu'})=
\begin{cases}
0&(v=1\text{ or }v=a+1)\\
1&(3\ls v\ls a-1).
\end{cases}
\]
\end{propn}

\begin{pf}
The construction of the homomorphism $\sigma$ above shows that the dimension of the homomorphism space is at least that claimed.  So we just have to show the reverse inequality. By Lemma \ref{815hom}, we have
\[
\dim_\bbf\hom_{\bbf\sss n}(S^\la,S^{\mu'})=\dim_\bbf\hom_{\bbf\sss n}(S^\mu,S^{\la'}),
\]
and we can use the technique outlined in \S\ref{homsec} to compute the right-hand side, since $\mu$ is $2$-regular. So suppose $\theta$ is a linear combination of semistandard homomorphisms $\hth T:S^\mu\to M^{\la'}$ such that $\ps dt\circ\theta=0$ for all $d,t$.

To begin with, we consider $\ps21\circ\hth T$ for each $T$.  Using Lemma \ref{lemma5}, this equals zero if $T$ has a $2$ in each row, because the homomorphisms occurring in Lemma \ref{lemma5} each appear with a coefficient $\binom21$, which is zero in $\bbf$. Otherwise, $\ps21\circ\hth T$ is either a single semistandard homomorphism or a sum of two semistandard homomorphisms.  Moreover, the semistandard tableaux that occur for the various $T$ are distinct.  Hence in order to have $\ps21\circ\theta=0$, $\theta$ can only involve semistandard homomorphisms $\hth T$ for those $T$ having a $2$ in each row. In particular, $\theta=0$ when $v=a+1$, since in this case there is only one semistandard tableau, whose first row consists entirely of $1$s.

Now we consider $\ps22\circ\hth T$ for each of these $T$. If $T$ has a $2$ and a $3$ in each row, we get $\ps22\circ\hth T=0$, while if $T$ has a $2$ in each row and two $3$s in the same row, $\ps22\circ\hth T$ is a semi\-standard homomorphism.  Again, all the semistandards that occur in this way are different, so $\theta$ cannot involve any tableaux of the latter type.  In particular, if $v=1$ then $\theta=0$.

Next consider $\ps d1\circ\hth T$ where $4\ls d<a$ and $T$ is a semistandard tableau having a $2$ and a $3$ in each row. $T$ contains a single $d$ and a single $d+1$.  If these both lie in the same row of $T$, then $\ps d1\circ\hth T=0$. Otherwise, $\ps d1\circ\hth T$ is a semistandard homomorphism $\hth{T'}$.  If $U$ is another semistandard tableau and $\ps d1\circ\hth U$ is a semistandard homomorphism $\hth{U'}$, then $T'=U'$ if and only if $U$ is obtained from $T$ by interchanging $d$ and $d+1$; hence any two such tableaux must occur in $\theta$ with equal coefficients. Applying this for all $d\gs4$ and all $T$, we find that $\theta$ must be a scalar multiple of the sum of all semistandard homomorphisms $\hth T$ for $T$ having a $2$ and a $3$ in each row. Hence $\dim_\bbf\hom_{\bbf\sss n}(S^\mu,S^\la)\ls1$, and we are done.
\end{pf}

\begin{eg}
We provide an example to illustrate the above proof for the benefit of the reader who may not be familiar with this technique.  We take $(a,b,u,v)=(4,6,8,5)$.  (In fact, the Specht module $S^{(8,5)}$ is reducible, so is ultimately irrelevant to our main theorem, but it serves well for this example.)

We suppose we have a linear combination $\theta$ of semistandard homomorphisms such that $\ps dt\circ\theta=0$ for all $d,t$. For this example, we abuse notation by identifying a tableau with the corresponding homomorphism. The first step of the proof is to eliminate most of the possible semistandard homomorphisms by taking $d=2$, $t=1$.  For example, Lemma \ref{lemma5} gives
\[
\ps21\circ\,\young(11113346,22578)
\,=\,
\young(11112346,22578),
\]
and no other semistandard tableau can give the semistandard tableau on the right in this way with non-zero coefficient; note that the tableau
\[
\young(11112346,23578)
\]
does give this tableau, but with a coefficient of $\binom21=0$. So since $\ps21\circ\theta=0$, our initial tableau cannot possibly occur in $\theta$.  Arguing in this way, one finds that the only semistandard tableaux which can occur in $\theta$ are those with a $2$ in each row, i.e.\ those of the form
\[
\young(1111233\star,2\star\star\star\star),\qquad
\young(111123\star\star,23\star\star\star)\quad\text{or}\quad
\young(11112\star\star\star,233\star\star).
\]
Now the first and last of these three types can be ruled out using the same argument with $\ps22$. So $\theta$ can only involve tableaux with a $2$ and a $3$ in each row; call these \emph{usable} tableaux.  Now we consider $\ps d1\circ\theta$ for $d\gs4$.  Now for each usable tableau $T$, $\ps d1\circ\hth T$ is either zero (if $d$ and $d+1$ occur in the same row in $T$) or a semistandard homomorphism.  Furthermore, these semistandard homomorphisms `pair up'; for example, with $d=4$ we have
\[
\ps41\circ\,\young(11112356,23478)\,=\ps41\circ\,
\young(11112346,23578)\,=\,
\young(11112346,23478).
\]
Since the semistandard tableau on the right can only arise in this way from the two semistandard tableaux on the left, these two semistandard homomorphisms must occur with equal coefficients in $\theta$.  Now we observe that we can get from any usable tableau to any other by a sequence of steps in which we interchange the integers $d,d+1$ for various values of $d\gs4$.  So if we apply the above argument for all $d\gs4$, we see that all usable tableaux occur with the same coefficient in $\theta$.
\end{eg}

\subsection{Homomorphisms from $S^\mu$ to $S^\la$}

Now we consider homomorphisms from $S^\mu$ to $S^\la$, where $\la,\mu$ are as above. In view of Proposition \ref{cdhomdim1}, \emph{we assume for the rest of this section that $3\ls v\ls a-1$}. It turns out that all such homomorphisms can be expressed as linear combinations of $\hth A$ and $\hth B$, where $A,B$ are the following $\mu$-tableaux of type $\la$:
\begin{align*}
A&=
\footy{\gyoungx(1.2,;1_2\hdss;1;2;2;2;3;4_2\hdss;\bpt,;1_2\hdss;1;1;1)};\\
B&=
\footy{\gyoungx(1.2,;1_2\hdss;1;1;1;2;3;4_2\hdss;\bpt,;1_2\hdss;1;2;2)}.
\end{align*}
Note that our assumptions on the parameters $a,b,u,v$ mean that these tableaux really do exist, i.e.\ there are enough $1$s to fill the bottom row.

\begin{lemma}\label{abnz}
$\hth A$ and $\hth B$ are non-zero, and are linearly independent if $v\ls b+1$.
\end{lemma}

\begin{pf}
It is straightforward to express $\hth A$ and $\hth B$ as linear combinations of semistandard homomorphisms using a single application of Lemma \ref{lemma7}; in each case we get at least one semistandard appearing, so the homomorphisms are non-zero.  If in addition $v\ls b+1$, then in the expression for $\hth A$ there is at least one semistandard tableau with two $2$s in the first row; there is no such tableau appearing in the expression for $\hth B$, so $\hth A,\hth B$ are linearly independent.
\end{pf}

\begin{propn}\label{abhoms}\leavevmode
\begin{itemize}
\item
\vspace{-\topsep}
If $a-v\equiv3\ppmod4$ or $v=b+3$, then $\ps dt\circ\hth A=0$ for all admissible $d,t$.
\item
If $v\equiv1\ppmod4$, then $\ps dt\circ\hth B=0$ for all admissible $d,t$.
\item
If $a\equiv0\ppmod4$, then $\ps dt\circ(\hth A+\hth B)=0$ for all admissible $d,t$.
\end{itemize}
\end{propn}

\begin{pf}
Lemma \ref{lemma5} immediately gives $\ps d1\circ\hth A=\ps d1\circ\hth B=0$ for $d\gs2$.  Using the fact that $A,B$ each have an odd number of $1$s in each row, we also get \[
\ps1t\circ\hth A=\ps1t\circ\hth B=0
\]
for $t=1,3$.  Finally, we have $\ps12\circ\hth B=0$ if $v\equiv1\ppmod 4$, and $\ps12\circ\hth A=0$ if $a-v\equiv3\ppmod4$ or $v=b+3$ (where we apply Lemma \ref{lemma7} in the latter case), and $\ps12\circ\hth A=\ps12\circ\hth B$ if $a\equiv0\ppmod4$.
\end{pf}

\begin{propn}\label{muladim}
\[
\dim_\bbf\hom_{\bbf\sss n}(S^\mu,S^{\la})=
\begin{cases}
2&(\text{if }a\equiv0\ppmod4,\ v\equiv1\ppmod4\text{ and }v\ls b+1)\\
1&(\text{otherwise}).
\end{cases}
\]
\end{propn}

\begin{pf}
By Lemma \ref{abnz} and Proposition \ref{abhoms} the dimension of the homomorphism space is at least that claimed.  Now we show the reverse inequality, by considering linear combinations of semistandard homomorphisms.

Throughout this proof, we'll write $\ca i$ for the set of semistandard $\mu$-tableaux of type $\la$ having exactly $i$ $2$s in the first row, for $i=0,1,2,3$, and let $\tau_i=\sum_{T\in\ca i}\hth T$.

Suppose we have a linear combination $\theta$ of semistandard homomorphisms $\hth T:S^\mu\to M^\la$ such that $\ps dt\circ\theta=0$ for all applicable $d,t$.

First consider $\ps d1\circ\hth T$ for $T\in\semi\mu\la$ and $d\gs3$. By Lemma \ref{lemma5}, $\ps d1\circ\hth T$ is either zero or a semistandard homomorphism (according to whether the $d$ and the $d+1$ in $T$ occur in the same row).  If it is non-zero, then there is exactly one other $T'\in\semi\mu\la$ such that $\ps d1\circ\hth T=\ps d1\circ\hth{T'}$, namely the tableau obtained by interchanging the $d$ and the $d+1$ in $T$. Hence $\hth T$ and $\hth{T'}$ must occur with the same coefficient in $\theta$.  Applying this for all $d\gs 3$, we find that for a fixed $i\in\{0,1,2,3\}$, all the homomorphisms $\hth T$ for $T\in\ca i$ occur with the same coefficient in $\theta$.  In other words, $\theta$ is a linear combination of $\tau_0,\tau_1,\tau_2,\tau_3$.

We can apply a similar argument in which we consider $\ps31\circ\hth T$ for $T\in\semi\mu\la$. Again $\ps31\circ\hth T$ is either zero or a semistandard homomorphism; and if it is non-zero, then the only other $T'$ having $\ps31\circ\hth{T'}=\ps31\circ\hth T$ is the tableau obtained by exchanging the $3$ in $T$ with a $2$ in the other row. $\hth T$ and $\hth{T'}$ occur with the same coefficient in $\theta$, and we deduce that $\theta$ must be a linear combination of $\tau_0+\tau_1$ and $\tau_2+\tau_3$.

Finally, we consider $\ps12\circ\theta$. Each $\mu$-tableau $T$ of type $(a+2,1^{b+1})$ contains a single $2$; let $\phi$ denote the sum of $\hth T$ for all those $T$ having the $2$ in row $1$, and $\chi$ the sum of all $\hth T$ for $T$ having the $2$ in row $2$. Using Lemma \ref{lemma5} and Lemma \ref{lemma7} (and recalling that $a$ is even and $v$ is odd), we have
\begin{align*}
\ps12\circ\tau_0&=\mbinom{v-1}2\chi,\\
\ps12\circ\tau_1&=\mbinom v2\phi,\\
\ps12\circ\tau_2&=\left(\mbinom{a+2}2+1\right)\phi+\chi,\\
\ps12\circ\tau_3&=\mbinom{a+2}2\phi.
\end{align*}
So if $v\equiv 3\ppmod4$, then $\ps12\circ(\tau_0+\tau_1)\neq0$, so $\theta$ cannot equal $\tau_0+\tau_1$. If $a\equiv2\ppmod4$, then $\ps12\circ(\tau_2+\tau_3)\neq0$, so $\theta$ cannot be $\tau_2+\tau_3$. Hence $\dim_\bbf\hom_{\bbf\sss n}(S^\mu,S^\la)\ls1$ in these cases. We also have $\dim_\bbf\hom_{\bbf\sss n}(S^\mu,S^\la)\ls1$ in the case where $b=d-3$, since in this case $\ca 2$ and $\ca 3$ are empty, so $\tau_2+\tau_3=0$.
\end{pf}

\subsection{Composing the homomorphisms}

Now we complete the analysis of when $S^\mu$ is a summand of $S^\la$, by composing the homomorphisms from the preceding subsections. This will be straightforward, using Proposition \ref{tabcomp}.

Recall that the space of homomorphisms from $S^\la$ to $S^{\mu'}$ is one-dimensional, spanned by the homomorphism $\sigma=\sum_{T\in\ct}\hth T$. On the other hand, the space of homomorphisms from $S^\mu$ to $S^\la$ has dimension one or two, each homomorphism being a linear combination of the homomorphisms $\hth A$ and $\hth B$. So it suffices to compute the compositions $\sigma\circ\hth A$ and $\sigma\circ\hth B$.

Let $D$ be the $\mu$-tableau
\[
D=\gyoung(;1;2;3_4\hdssss;u,;1;2;3_2\hdss;v)
\]
of type $\mu'$.  Then we have the following.

\begin{lemma}\label{uab}
Suppose $T\in\ct$, and let $x$ be the entry in the $(2,2)$-position of $T$. Then
\[
\hth T\circ\hth A=\hth D,
\qquad
\hth T\circ\hth B=\begin{cases}
\hth D&(x\ls v)\\
0&(x>v).
\end{cases}
\]
Furthermore, $\hth D\neq0$.
\end{lemma}

\begin{pf}
The fact that $\hth D\neq0$ is a simple application of Lemma \ref{lemma7}. To show that the compositions of homomorphisms are as claimed, take $T\in\ct$ and recall the notation of Proposition \ref{tabcomp}, with $S$ equal to either $A$ or $B$.

Suppose $X\in\calx$. Since each $T^i$ is a proper set, each $X^{ij}$ must be as well.  This means that if some integer $i$ appears in two sets $X^{kj},X^{lj}$, then the multinomial coefficient $\mbinom{X^{1j}_j+X^{2j}_j+X^{3j}_j+\dots}{X^{1j}_j,X^{2j}_j,X^{3j}_j,\dots}$ from Proposition \ref{tabcomp} will include a factor $\binom21$, which gives $0$.

So in order to get a non-zero coefficient in Proposition \ref{tabcomp}, we must have $X^{1j},X^{2j},X^{3j},\dots$ pairwise disjoint for each $j$, which means that we will have
\[
X^{11}\sqcup X^{21}\sqcup\dots=\{1,\dots,u\},\qquad X^{12}\sqcup X^{22}=\{1,\dots,v\};\tag*{($\dagger$)}
\]
so $U_X$ will equal $D$.

If $S=A$, the only way to achieve this is to have
\[
X^{11}=T^1\setminus\{1,\dots,v\},\quad X^{12}=\{1,\dots,v\},\quad X^{i1}=T^i\text{ for }i\gs2.
\]
Thus we have $\hth T\circ\hth A=\hth D$.

In the case $S=B$, let $y$ be the $(2,3)$-entry of $T$.  Then $y>x$. $X^{22}$ must contain either $x$ or $y$, so if $x>v$ then we cannot possibly achieve ($\dagger$). So we get $\hth T\circ\hth B=0$ in this case. If $x\ls v<y$, then the only way to achieve ($\dagger$) is to have $X^{22}=\{1,x\}$ and $X^{12}=\{2,\dots,\hat x,\dots,v\}$, and this yields $\hth T\circ\hth B=\hth D$.  Finally, if $y\ls v$, then there are three possible ways to achieve ($\dagger$); each of these gives a coefficient of $1$, and again we have $\hth T\circ\hth B=\hth D$.
\end{pf}

This result is very helpful: it tells us that the composition of $\sigma$ with a combination of $\hth A$ and $\hth B$ is a scalar multiple of $\hth D$; hence this composition is non-zero if and only if this scalar is non-zero. In order to use this result, we need to find the number of tableaux in $\ct$, and also the number of tableaux in $\ct$ in which the $(2,2)$-entry is at most $v$.  This is a straightforward count.

\Needspace*{2\baselineskip}\begin{lemma}\label{countu}\indent
\begin{itemize}
\item
\vspace{-\topsep}
The number of tableaux in $\ct$ is $\mbinom{u-v}{a-v}\mbinom{u+v-a-1}2$.
\item
The number of tableaux in $\ct$ whose $(2,2)$-entry is greater than $v$ is $\mbinom{u-v}{a-v}\mbinom{u-a}2$.
\end{itemize}
\end{lemma}

\begin{pf}[Proof of Theorem \ref{main}(1)]
Suppose $S^\mu=S^{(u,v)}$ is irreducible, with $u+v=a+b+3$.  Throughout this proof, all congruences are modulo $4$.  

Suppose first that $u,v$ satisfy the given conditions, i.e.\ $v\equiv3$ and $\mbinom{u-v}{a-v}$ is odd. The second condition implies in particular that $0\ls a-v\ls u-v$, which gives $v\ls\min\{a-1,b+3\}$; so the assumptions in force in this section are satisfied. In addition, Theorem \ref{irrspecht} gives $u\equiv2$.

We need to show that there are homomorphisms $S^\mu\stackrel\gamma\longrightarrow S^\la\stackrel\delta\longrightarrow S^{\mu'}$ such that $\delta\circ\gamma\neq0$. Since $3\ls v\ls a-1$, we can take $\delta=\sigma$.

If $a\equiv0$, take $\gamma=\hth A+\hth B$. By Proposition \ref{abhoms}, $\gamma$ is a homomorphism from $S^\mu$ to $S^\la$. By Lemma \ref{uab} and Lemma \ref{countu},
\[
\delta\circ\gamma=\mbinom{u-v}{a-v}\mbinom{u-a}2\hth D.
\]
The first term is odd by assumption; the second term is odd because $u-a\equiv2$, and $\hth D\neq0$ by Lemma \ref{uab}.

If $a\equiv2$, take $\gamma=\hth A$. Then $\gamma$ is a homomorphism from $S^\mu$ to $S^\la$, and
\[
\delta\circ\gamma=\mbinom{u-v}{a-v}\mbinom{u+v-a-1}2\hth D.
\]
Again, the first term is odd by assumption, the second term is odd because now $u+v-a-1\equiv2$, and $\hth D\neq0$.

Conversely, suppose we have homomorphisms $\gamma,\delta$ such that $\delta\circ\gamma\neq0$. By Proposition \ref{cdhomdim1} we can assume that $3\ls v\ls a-1$ and take $\delta=\sigma$.  From Proposition \ref{muladim} we can take $\gamma$ to be $\hth A$, $\hth B$ or $\hth A+\hth B$, according to the congruences in Proposition \ref{abhoms}.  Then $\delta\circ\gamma$ will be a scalar multiple of $\hth D$, and the scalar will include $\mbinom{u-v}{u-a}$ as a factor.  So this binomial coefficient must be odd, and all that remains is to show that $v\equiv3\ppmod4$.  We consider the three cases of Proposition \ref{abhoms}. Note that because $v>1$, Theorem \ref{irrspecht} gives $u-v\equiv3$.
\begin{description}
\item[\fbox{\rm $a-v\equiv3$ or $v=b+3$, $\gamma=\hth A$}]
In this case the coefficient of $\hth D$ in $\delta\circ\gamma$ is
\[
\mbinom{u-v}{a-v}\mbinom{u+v-a-1}2.
\]
The second binomial coefficient must be odd, so $u+v-a\equiv3$.  In the case $a-v\equiv3$, this is the same as saying $u\equiv2$, so that $v\equiv3$. In the case $v=b+3$, we have $a=u$, so that again $v\equiv3$.
\item[\fbox{\rm $a\equiv0$, $\gamma=\hth A+\hth B$}]
Now the coefficient of $\hth D$ in $\delta\circ\gamma$ is
\[
\mbinom{u-v}{a-v}\mbinom{u-a}2.
\]
The second binomial coefficient is odd only if $u\equiv a+2$, which is the same as saying $v\equiv 3$.
\item[\fbox{\rm $v\equiv1$, $\gamma=\hth B$}]
In this case the coefficient of $\hth D$ in $\delta\circ\gamma$ is
\[
\mbinom{u-v}{a-v}\left(\mbinom{u+v-a-1}2+\mbinom{u-a}2\right).
\]
Since $v\equiv1$, we have $u+v-a-1\equiv u-a$, so that $\mbinom{u+v-a-1}2$ and $\mbinom{u-a}2$ have the same parity. Hence $\delta\circ\gamma=0$, a contradiction.\qedhere
\end{description}
\end{pf}

\section{Irreducible summands of the form $S^{(u,v,2)}$}\label{uv2sec}

In this section, we find when one of our Specht modules $S^{(a,3,1^b)}$ has a summand isomorphic to an irreducible Specht module of the form $S^{(u,v,2)}$, where $u$ is even and $v$ is odd.  By Theorem \ref{jreg} and Lemma \ref{reg}, $D^{(u,v,2)}$ cannot appear as a composition factor of $S^{(a,3,1^b)}$ unless $(u,v,2)\dom(a,3,1^b)\reg$.  So we may assume that this is the case, which is the same as saying $v\ls\min\{a-1,b+1\}$.  We set out notation and assumptions for this section.

\begin{mdframed}[innerleftmargin=3pt,innerrightmargin=3pt,innertopmargin=3pt,innerbottommargin=3pt,roundcorner=5pt,innermargin=-3pt,outermargin=-3pt]
\noindent\textbf{Assumptions and notation in force throughout Section \ref{uv2sec}:}

$\la=(a,3,1^b)$ and $\mu=(u,v,2)$, where $a,b,u,v$ are positive integers with $a,b,u$ even, $a\gs4$, $u>v>2$, $n=a+b+3=u+v+2$ and $v\ls\min\{a-1,b+1\}$.
\end{mdframed}

\subsection{Homomorphisms from $S^\la$ to $S^{\mu'}$}

We begin by constructing a homomorphism from $S^\la$ to $S^{\mu'}$.  As in \S\ref{hlamu'1}, we construct this using non-semistandard tableaux.

Let $\ct$ be the set of $\la$-tableaux having the form
\[
\gyoung(;1;2;3_\sq\hdsq;v;\star_\sq\hdsq;\star,;1;1;2,;2,;3,|\sq\vdsq,;v,;\star,|\sq\vdsq,;\star)
\]
in which the $\star$s represent the numbers from $v+1$ to $u$, and in which the entries are weakly increasing along the first row and down the first column.  Let $\sigma=\sum_{T\in\ct}\hth T$.

\begin{propn}\label{sigmahom2}
With the notation and assumptions above, we have $\ps dt\circ\sigma=0$ for all $d,t$.
\end{propn}

\begin{pf}
For $d\gs v$ and $t=1$, we use the same argument as that used in several proofs above: for $T\in\ct$ either $\ps d1\circ\hth T=0$, or there is a unique other $T'\in\ct$ with $\ps d1\circ\hth{T'}=\ps d1\circ\hth T$.

The cases where $2\ls d\ls v$ are easier: in this case Lemma \ref{lemma5} and Lemma \ref{lemma7} imply that we have $\ps dt\circ\hth T=0$ for all $T\in\ct$.

So we are left with the cases where $d=1$ and $t\in\{1,2,3\}$.  For $T\in \ct$ we have $\ps 13\circ\hth T$ immediately from Lemma \ref{lemma5}, while $\ps12\circ\hth T$ is a homomorphism labelled by a tableau of the form
\[
\gyoung(;1;2;3_\sq\hdsq;v;\star_\sq\hdsq;\star,;1;1;1,;1,;3,|\sq\vdsq,;v,;\star,|\sq\vdsq,;\star).
\]
But this homomorphism is zero by Lemma \ref{lemma7}.  Finally, $\ps11\circ\hth T$ is the sum of the homomorphisms labelled by two tableaux
\[
\gyoung(;1;2;3_\sq\hdsq;v;\star_\sq\hdsq;\star,;1;1;1,;2,;3,|\sq\vdsq,;v,;\star,|\sq\vdsq,;\star),\qquad
\gyoung(;1;2;3_\sq\hdsq;v;\star_\sq\hdsq;\star,;1;1;2,;1,;3,|\sq\vdsq,;v,;\star,|\sq\vdsq,;\star).
\]
But these two homomorphisms are equal by Lemma \ref{lemma7}, and we are done.
\end{pf}

Now, as in Section \ref{hlamu'1} we have to show that $\sigma\neq0$.  Again, we use a dominance argument.

\begin{propn}\label{sigmanz2}
With the notation above, $\sigma\neq0$.
\end{propn}

\begin{pf}
We'll show that when $\sigma$ is expressed as a linear combination of semistandard homomorphisms, the homomorphism $\hth S$ occurs with non-zero coefficient, where
\[
S=
\footy{\gyoungx(1.2,;1;1;1;2;4;5_\sq\hdsq;v;\bph_\sq\hdsq;u,;2;2;3,;3,|2\vdss,;\bpt)}.
\]
Recall that when $\hth T$ is expressed as a linear combination of semistandard homomorphisms, the coefficient of $\hth S$ is zero unless $S\dom T$.  The only elements of $\ct$ which are dominated by $S$ are the tableaux of the form
\[
T[i]=
\footy{\gyoungx(1.2,;1;2;3;4_\sq\hdsq;v;i;\bph_\sq\hdsq;u,;1;1;2,;2,;3,;\vds,;v,;\vpo,;\vds,;\hati,;\vds,;\bpt)}
\]
for $v+1\ls i\ls b+2$. Consider applying Lemma \ref{lemma7} to $T[i]$, to move the two $1$s from row $2$ to row $1$. Of the tableaux appearing in that lemma with non-zero coefficient, the only ones dominated by $S$ are those having no more than four entries less than $4$ in the first row; these are the tableaux $T'[i]$ and $T'[i,j]$ for $4\ls j\ls v$, where
\begin{align*}
T'[i]&=
\footy{\gyoungx(1.2,;1;1;1;2;4_\sssq\hdsssq;v;\bph_\sq\hdsq;u,;2;3;i,;2,;3,;\vds,;v,;\vpo,;\vds,;\hati,;\vds,;\bpt)},
\\
T'[i,j]&=
\footy{\gyoungx(1.2,;1;1;1;2;4;\hds;\hatj;\hds;v;i;\bph;\hds;u,;2;3;j,;2,;3,;\vds,;v,;\vpo,;\vds,;\hati,;\vds,;\bpt)}.
\end{align*}
So, modulo homomorphisms labelled by tableaux not dominated by $S$, we have $\sigma=\sum_i\hth{T'[i]}+\sum_{i,j}\hth{T'[i,j]}$.  However, two applications of Lemma \ref{lemma7} show that $\hth{T'[i,j]}=0$ for all $i,j$, and Lemma \ref{lemma7} also gives $\hth{T'[i]}=\hth S$.

So the coefficient of $\hth S$ in $\sigma$ is $b+2-v$, which is odd; so $\sigma\neq0$.
\end{pf}

As before, we find that $\sigma$ is the only homomorphism from $S^\la$ to $S^{\mu'}$ up to scaling.

\begin{propn}\label{cd2homdim1}
With $\la,\mu$ as above,
\[
\dim_\bbf\hom_{\bbf\sss n}(S^\la,S^{\mu'})=1.
\]
\end{propn}

\begin{pf}
The existence of the homomorphism $\sigma$ shows that the space of homomorphisms is non-zero. To show that it has dimension at most $1$, we again use the fact that
\[
\dim_\bbf\hom_{\bbf\sss n}(S^\la,S^{\mu'})=\dim_\bbf\hom_{\bbf\sss n}(S^\mu,S^{\la'}).
\]
So suppose $\theta$ is a linear combination of semistandard homomorphisms $\hth T:S^\mu\to S^{\la'}$ such that $\ps dt\circ\theta=0$ for all $d,t$.

First of all, consider $\ps21\circ\theta$.  Given a semistandard tableau $T$, we can use Lemma \ref{lemma5} to compute $\ps21\circ\hth T$, and then if necessary use Lemma \ref{lemma7} (to move a $2$ from row $3$ to row $2$) to express this composition as a linear combination of semistandard homomorphisms. We find that if $T$ has two $2$s in its first row, then $\ps21\circ\hth T$ involves a semistandard tableau which does not occur in any other $\ps21\circ\hth{T'}$; hence the coefficient of $\hth T$ in $\theta$ must be zero.

Now we do the same thing with $\ps22$: in this case we find that if $T$ is a semistandard tableau having two $3$s in its first row, then $\ps22\circ\hth T$ involves a semistandard homomorphism which does not occur in any other $\ps22\circ\hth{T'}$ (except possibly for a tableau $T'$ already ruled out in the paragraph above).  So we may restrict attention to those $T$ having at most one $2$ and one $3$ in the first row.

Now return to $\ps21\circ\hth T$, for $T$ of the form
\begin{align*}
&\gyoung(;1_5\hdsssss;1;2;3;\xone_2\hdss;\xs,;2;\zone;\ztwo;\hds;\zt,;3;k),
\\
\intertext{where $x_1,\dots,x_s,z_1,\dots,z_t,k$ are the integers $4,\dots,a$ in some order. When we express $\ps21\circ\hth T$ as a linear combination of semistandard homomorphisms, we find that the homomorphism labelled by}
&\gyoung(;1_5\hdsssss;1;2;3;\xone_2\hdss;\xs,;2;2;\ztwo;\hds;\zt,;\zone;k)
\end{align*}
occurs with non-zero coefficient; but this homomorphism does not occur in any other $\ps21\circ\hth{T'}$ (except for $T'$ having two $3$s in its first row). So for any $T$ of the above form, the coefficient of $\hth T$ in $\theta$ must be zero.

Now any semistandard tableau which contributes to $\theta$ must be of one of the following eight forms.
\begin{alignat*}2
1.\quad&
\gyoung(;1;1;1_5\hdsssss;1;2;\star_2\hdss;\star,;2;3;\star_3\hdsss;\star,;3;\star)
&\qquad
2.\quad&
\gyoung(;1;1;1_5\hdsssss;1;\star_3\hdsss;\star,;2;2;3;\star_2\hdss;\star,;3;\star)
\\
3.\quad&
\gyoung(;1;1;1_5\hdsssss;1;\star_3\hdsss;\star,;2;2;\star_3\hdsss;\star,;3;3)
&\qquad
4.\quad&
\gyoung(;1;1;1_5\hdsssss;1;3;\star_2\hdss;\star,;2;2;\star_3\hdsss;\star,;3;\star)
\\
5.\quad&
\gyoung(;1;1;1_5\hdsssss;1;3;\star_2\hdss;\star,;2;2;3;\star_2\hdss;\star,;\star;\star)
&\qquad
6.\quad&
\gyoung(;1;1;1_5\hdsssss;1;\star_3\hdsss;\star,;2;2;3;3;\star;\hds;\star,;\star;\star)
\\
7.\quad&
\gyoung(;1;1;1_5\hdsssss;1;2;3;\star;\hds;\star,;2;3;\star_3\hdsss;\star,;\star;\star)
&\qquad
8.\quad&
\gyoung(;1;1;1_5\hdsssss;1;2;\star_2\hdss;\star,;2;3;3;\star_2\hdss;\star,;3;\star)
\end{alignat*}
In each case, the $\star$s represent the numbers from $4$ to $a$.  Note that in each of these tableaux, the entries $4,\dots,a$ must all occur in different columns (the assumption $v\ls b+1$ means that any column of length at least two has a $1$ at the top).  So we can consider the homomorphisms $\ps d1$ for $d\gs3$, and repeat the argument used in the last paragraph of the proof of Proposition \ref{cdhomdim1}, to show that if $T,T'$ are two tableaux which have their $1$s, $2$s and $3$s in the same positions, then $\hth T$ and $\hth{T'}$ occur with the same coefficient in $\theta$.  Hence $\theta$ is a linear combination of the homomorphisms $\tau_1,\dots,\tau_8$, where $\tau_i$ is the sum of all homomorphisms $\hth T$ for $T$ of type $i$.

Once more we can consider $\ps21\circ\theta$: when we compute $\ps21\circ\tau_5$, we obtain (in addition to some other semistandard tableaux) the sum of the semistandard tableaux of the form
\[
\gyoung(;1;1;1_5\hdsssss;1;3;\star_2\hdss;\star,;2;2;2;\star_2\hdss;\star,;\star;\star),
\]
which do not occur in any other $\ps21\circ\tau_i$ (note these tableaux do occur when we compute $\ps21\circ\hth T$ for $T$ of type $4$, but each one occurs twice when we sum over tableaux of type $4$, so the contributions cancel).  So $\tau_5$ does not appear in $\theta$.

Next we consider $\ps31\circ\tau_i$ for each $i$.  For $i=3,6$ or $8$, we find that $\ps31\circ\tau_i$ involves semistandard tableaux which do not occur in any other $\ps31\circ\tau_i$; so $\tau_3,\tau_6,\tau_8$ cannot occur in $\theta$.  Moreover, we find that $\ps31\circ\tau_1=\ps31\circ\tau_7$ and $\ps31\circ\tau_2=\ps31\circ\tau_4$, and that these two homomorphisms are linearly independent.  So $\theta$ must be a linear combination of $\tau_1+\tau_7$ and $\tau_2+\tau_4$.

Finally we return once more to $\ps21\circ\theta$.  We find that $\ps21\circ\tau_1=\ps21\circ\tau_4\neq0$, while $\ps21\circ\tau_2=\ps21\circ\tau_7=0$.  So $\tau_1$ and $\tau_4$ must appear with the same coefficient in $\theta$; so $\theta$ must be a scalar multiple of $\tau_1+\tau_2+\tau_4+\tau_7$, and so the homomorphism space has dimension at most $1$.
\end{pf}

\subsection{Homomorphisms from $S^\mu$ to $S^\la$}

Now we consider homomorphisms from $S^\mu$ to $S^\la$.  We begin by constructing a non-zero homomorphism.  Let $C$ be the $\mu$-tableau
\[
\footy{\gyoungx(1.2,;1;1_\sssq\hdsssq;1;2;3_2\hdss;\bpt,;1;1_\sq\hdsq;1,;2;2)}
\]
of type $\la$.

\begin{propn}\label{cd2chom}
With $C$ as above, we have $\ps dt\circ\hth C=0$ for all $d,t$, and $\hth C\neq0$.
\end{propn}

\begin{pf}
Showing the first statement is very easy, using Lemma \ref{lemma5}.  The only homomorphisms that occur in that lemma with non-zero coefficient are labelled by tableaux with more than $v$ $1$s in rows $2$ and $3$, and therefore by Lemma \ref{lemma7} are zero.

Showing that $\hth C\neq0$ is also straightforward using Lemma \ref{lemma7}.  We apply this lemma to move the $1$s from row $2$ to row $1$, and then again to move the $2$s from row $3$ to row $2$.  The tableau
\[
\footy{\gyoungx(1.2,;1;1;1_\ssssq\hdssssq;1;\vpt_\sq\hdsq;\bpt,;2;2;2;5_\sq\hdsq;\vpo,;3;4)}
\]
(for example) labels a homomorphism occurring with non-zero coefficient.
\end{pf}

\begin{propn}\label{cd2homdim2}
With $\la,\mu$ as above,
\[
\dim_\bbf\hom_{\bbf\sss n}(S^\mu,S^\la)=1.
\]
\end{propn}

\begin{pf}
The existence of the homomorphism $\hth C$ above shows that the space of homomorphisms is non-zero.  So we just need to show the upper bound on the dimension. So suppose $\theta$ is a linear combination of semistandard homomorphisms $\hth T:S^\mu\to S^\la$ such that $\ps dt\circ\theta=0$ for all $d,t$.

For $3\ls d\ls b+1$, say that a semistandard tableau $T$ is \emph{$d$-bad} if the entries $d,d+1$ appear in the same column of $T$.  Note that this must be column $1$ or $2$, because the assumption $v\ls a-1$ guarantees that any column of length greater than $1$ has a $1$ at the top.

We claim that if $T$ is $d$-bad, then $\hth T$ cannot appear in $\theta$. To show this, we consider $\ps d1\circ\hth T$ for every semistandard $T$.  If $T$ is not $d$-bad, then by Lemma \ref{lemma5} $\ps d1\circ\hth T$ is either zero or a homomorphism labelled by a semistandard tableau with two $d$s in different rows.  If $T$ is $d$-bad, then we can express $\ps d1\circ\hth T$ as a linear combination of semistandard homomorphisms using Lemma \ref{lemma5} together with Lemma \ref{newsemilem}. For example, if
\begin{align*}
T&=\young(11111122248\ten,369\eleven\twelve,57)
\\
\intertext{then $\ps61\circ\hth T=\hth{T'}$, where}
T'&=\young(11111122248\ten,369\eleven\twelve,56)
\end{align*}
and we can semistandardise this using Lemma \ref{newsemilem}, taking
\[
A=\{3\},\qquad B=\{5,6,7,9,11,12\},\qquad C=\emptyset
\]
to express $\hth{T'}$ as a sum of fourteen semistandard homomorphisms.

Doing this for each $d$-bad tableau $T$, we find that $\ps d1\circ\hth T$ is a sum of homomorphisms labelled by semistandard tableaux with the same first row as $T$; furthermore, at least one of these tableaux will have two $d$s in the second row.  Moreover, each $d$-bad tableau will yield a tableau of this kind which does not occur for any other $d$-bad tableau $T'$.  To see this, suppose first of all that $d,d+1$ occur in the second column of $T$.  Then there is no other $d$-bad tableau with the same first row as $T$, so any tableau occurring in $\ps d1\circ\hth T$ with two $d$s in the second row can only possibly occur in $\ps d1\circ\hth T$.  Alternatively, if $d,d+1$ occur in the first column of $T$, then $T$ has the form
\[
\footy{\gyoungx(1.2,;1_5\hdsssss;1;2;2;2;\xone_\sq\hdsq;\xs,;d;\zone;\ztwo;\hds;\zt,;\dpo;k)}.
\]
There are $v-2$ other $d$-bad tableaux with the same first row as $T$, and they all also have the same $(2,2)$-entry as $T$.  Hence when we apply Lemma \ref{lemma5} and Lemma \ref{newsemilem} (or equivalently Lemma \ref{lemma7}), we find that the homomorphism labelled by
\[
\gyoung(;1_5\hdsssss;1;2;2;2;\xone_\sq\hdsq;\xs,;d;d;\ztwo;\hds;\zt,;\zone;k)
\]
occurs in $\ps d1\circ\hth T$ but not in $\ps d1\circ\hth{T'}$ for any other $T'$.

So in $\theta$ the coefficient of $\hth T$ is zero for any $d$-bad tableau.  In particular, this means that for any $\hth T$ occurring in $\theta$, $\ps d1\circ\hth T$ is either zero or a single semistandard homomorphism. Now we claim that the coefficient of $\hth T$ is zero whenever $T$ has two or three $2$s in its first row.  Supposing this is false, take a $T$ with at least two $2$s in its first row such that $\hth T$ appears with non-zero coefficient in $\theta$, and suppose that $T$ is minimal (with respect to the dominance order) subject to this property.  Suppose the $(3,1)$-entry of $T$ is greater than $3$; then this entry equals $d+1$ for some $d\gs3$, and the entry equal to $d$ cannot be the $(2,1)$-entry (because $T$ is not $d$-bad).  So the $d$ and the $d+1$ in $T$ lie in different rows and different columns, and $\ps d1\circ\hth T$ is the semistandard homomorphism obtained by replacing the $d+1$ with a $d$.  The only other semistandard tableau $T'$ such that $\ps d1\circ\hth{T'}=\ps d1\circ\hth 
T$ is the tableau obtained by interchanging the $d$ and the $d+1$ in $T$, so $\hth{T'}$ must also occur with non-zero coefficient.  But $T\dom T'$, contradicting the choice of $T$.

So the $(3,1)$-entry in $T$ must be $3$ (and hence the $(2,1)$-entry is $2$).  Now consider the $(3,2)$-entry; call this $d+1$.  Then $d\gs4$, and the $d$ in $T$ cannot occur in the $(2,2)$-position (because $T$ is not $d$-bad).  So we can repeat the above argument and show that that there is a tableau $T'\domby T$ such that $\hth{T'}$ occurs in $\theta$; contradiction.

We now know that every semistandard homomorphism occurring in $\theta$ has at least two $2$s in the second row.  This means in particular that the entries $3,\dots,b+2$ lie in different columns.  So we can repeat the argument from Proposition \ref{cdhomdim1} and Proposition \ref{cd2homdim1} to show that $\theta$ must be a linear combination of $\tau_0$ and $\tau_1$, where $\tau_i$ is the sum of all homomorphisms labelled by semistandard tableau with $i$ $2$s in the first row.  If $u=a$, then $\tau_1=0$, and so the space of homomorphisms $S^\mu\to S^\la$  has dimension at most $1$.  if $u>a$, then $\ps21\circ\tau_0=\ps21\circ\tau_1\neq0$, so $\theta$ must be a scalar multiple of $\tau_0+\tau_1$ and again the homomorphism space has dimension at most $1$.
\end{pf}

\subsection{Composing the homomorphisms}

We have constructed homomorphisms $S^\mu\stackrel{\hth C}\longrightarrow S^\la\stackrel\sigma\longrightarrow S^{\mu'}$, and shown that these homomorphisms are unique up to scaling. To complete this section, we just need to compute the composition of these homomorphisms.

Let $E$ be the $\mu$-tableau
\[
\footy{\gyoungx(1.2,;1;2;3_\sq\hdsq;v;\vpo_\sq\hdsq;u,;1;2;3_\sq\hdsq;v,;1;2)}
\]
of type $\mu'$. Then we have the following.

\begin{propn}\label{comp2}
For $T\in\ct$, we have $\hth T\circ\hth C=\hth E\neq0$, and therefore we have $\sigma\circ\hth C\neq0$ if and only if $\mbinom{u-v}{a-v}$ is odd.
\end{propn}

\begin{pf}
It is easy to express $\hth E$ as a linear combination of semistandard homomorphisms using three applications of Lemma \ref{lemma7}, from which it follows that $\hth E\neq0$.

To prove that $\hth T\circ\hth C=\hth E$, use the notation of Proposition \ref{tabcomp}, with $S=C$.  Suppose $X\in\calx$ is such that the coefficient of $\hth{U_X}$ in Proposition \ref{tabcomp} is non-zero.  Since $X^{31}$ must be $\{2\}$, we cannot have $X^{21}=\{2\}$ (because this would give a factor $\binom21$), so $X^{21}=\{1\}$ and hence $X^{23}=\{1,2\}$.  Now if $X^{11}$ contains any of the numbers $1,\dots,v$ then again we get a factor $\binom21$.  So we have $X^{12}=\{1,\dots,v\}$, which determines $X$, and we find that $\hth T\circ\hth C=\hth E$ as required.

So we have $\sigma\circ\hth C=|\ct|\hth E$, which is non-zero if and only if $|\ct|$ is odd.  But it is easy to see that $|\ct|=\mbinom{u-v}{a-v}$, and the proposition is proved.
\end{pf}

\begin{pf}[Proof of Theorem \ref{main}(2)]
Suppose $S^\mu=S^{(u,v,2)}$ is irreducible.

Suppose first that $\mbinom{u-v}{a-v}$ is odd. This implies in particular that $0\ls a-v\ls u-v$, so $v\ls\min\{a-1,b+1\}$. So the assumptions of this section are valid, and we have homomorphisms $\hth C:S^\mu\to S^\la$ and $\sigma:S^\la\to S^{\mu'}$. By Proposition \ref{comp2}, $\sigma\circ\hth C\neq0$, so $S^\mu$ is a summand of $S^\la$.

Conversely, suppose we have homomorphisms $S^\mu\stackrel\gamma\longrightarrow S^\la\stackrel\delta\longrightarrow S^{\mu'}$ with $\delta\circ\gamma\neq0$. By Propositions \ref{cd2homdim1} and \ref{cd2homdim2}, $\delta$ must be a scalar multiple of $\sigma$, and $\gamma$ must be a scalar multiple of $\hth C$. Hence by Proposition \ref{comp2}, $\mbinom{u-v}{a-v}$ is odd.
\end{pf}

\section{Decomposability of Specht modules}\label{whichdec}

In this section, we prove Corollary \ref{maincor}, which answers the question of which Specht modules are shown to be decomposable by Theorem \ref{main}. First we consider the case where $a+b\equiv0\ppmod8$.

\begin{propn}\label{ab0}
Suppose $n\equiv3\ppmod8$, and $a$ is even, with $6\ls a\ls n-7$.  Let $b=n-a-3$.  Then $S^{(a,3,1^b)}$ has an irreducible summand of the form $S^{(u,v,2)}$.
\end{propn}

\begin{pf}
Using Theorem \ref{main}, we need to show that there is a pair $u,v$ with $u+v+2=n$ such that $S^{(u,v,2)}$ is irreducible and $\mbinom{u-v}{u-a}$ is odd.  By Theorem \ref{irrspecht}, $(u,v,2)$ is irreducible if and only if
\[
v\equiv1\pmod4,\qquad u-v\equiv-1\pmod{2^{l(v-2)}},
\]
where $l(k)=\lceil\log_2(k+1)\rceil$ for an integer $m$.

We use induction on $n$, with our main tool being the following well-known relations modulo $2$ on binomial coefficients:
\[
\binom{2x}{2y}\equiv\binom{2x+1}{2y}\equiv\binom{2x+1}{2y+1}\equiv\binom xy,\qquad \binom{2x}{2y+1}\equiv0\pmod 2.
\]
We consider three cases.
\begin{description}
\item[\fbox{$a=6,8,n-9$ or $n-7$}]
In this case, take $v=5$ (so $u=n-7$).  Since $n\equiv3\ppmod4$, we get $u\equiv0\ppmod4$, which means that $u-v\equiv3\ppmod4$, so $S^{(u,v,2)}$ is irreducible. Furthermore, the binomial coefficients
\[
\binom{u-5}0,\binom{u-5}1,\binom{u-5}2,\binom{u-5}3
\]
are all odd, which means that $\mbinom{u-5}{u-a}$ must be odd.
\item[\fbox{$n\equiv11\ppmod{16}$}]
In this case, let
\[
n'=\frac{n+11}2,\quad a'=
\begin{cases}
\mfrac{a+6}2&(a\equiv2\ppmod4)\\[5pt]
\mfrac{a+4}2&(a\equiv0\ppmod4).
\end{cases}
\]
Then $n',a'$ satisfy the conditions of the proposition, and $n'<n$ (note that the conditions on $a$ mean that $n>11$). So by induction there is a pair $u',v'$ such that
\[
v'\equiv1\pmod4,\qquad u-v\equiv-1\pmod{2^{l(v'-2)}},\qquad\mbinom{u'-v'}{u'-a'}\equiv1\pmod2.
\]
Note that because $u'-v'$ is odd and $u'-a'$ is even, this also gives $\mbinom{u'-v'}{u'-a'+1}$ odd.

We let $u=2u'-4$ and $v=2v'-5$.  Then $u+v+2=n$, and we have $v\equiv1\ppmod4$ and
\[
u-v=2(u'-v')+1\equiv-1\ppmod{2^{l(v'-2)+1}},
\]
with $l(v-2)\ls l(v'-2)+1$.  So $S^{(u,v,2)}$ is irreducible.  Furthermore
\[
\binom{u-v}{u-a}=\binom{2u'-2v'+1}{2u'-2a'(+2)}\equiv\binom{u'-v'}{u'-a'(+1)}\equiv1\pmod2,
\]
and we are done.
\item[\fbox{$n\equiv3\ppmod{16}$, $10\ls a\ls n-11$}]
In this case, let
\settowidth\frt{$\mfrac{a+2}2$}
\[
n'=\frac{n+3}2,\quad a'=
\begin{cases}
\mfrac{a+2}2&(a\equiv2\ppmod4)\\[5pt]
\hbox to \frt{\hfil$\mfrac a2$\hfil}&(a\equiv0\ppmod4).
\end{cases}
\]
Then $n',a'$ satisfy the conditions of the proposition, and $n'<n$. So by induction there is a pair $u',v'$ such that
\[
v'\equiv1\pmod4,\qquad u-v\equiv-1\pmod{2^{l(v'-2)}},\qquad\mbinom{u'-v'}{u'-a'}\equiv1\pmod2.
\]
Again, because $u'-v'$ is odd and $u'-a'$ is even, $\mbinom{u'-v'}{u'-a'+1}$ is also odd.

We let $u=2u'$ and $v=2v'-1$.  Then $u+v+2=n$, $v\equiv1\ppmod4$, and
\[
u-v=2(u'-v')+1\equiv-1\pmod{2^{l(v'-2)+1}},
\]
with $l(v-2)\ls l(v'-2)+1$.  So $S^{(u,v,2)}$ is irreducible.  Furthermore
\[
\binom{u-v}{u-a}=\binom{2u'-2v'+1}{2u'-2a'(+2)}\equiv\binom{u'-v'}{u'-a'(+1)}\equiv1\pmod2.\tag*{\raisebox{-10pt}{\qedhere}}
\]
\end{description}
\end{pf}

The next result addresses most of the cases where $a+b\equiv2\ppmod8$.

\begin{propn}\label{ab1}
Suppose $n\equiv5\ppmod8$, and $a$ is even, with $8\ls a\ls n-7$.  Let $b=n-a-3$.  Then $S^{(a,3,1^b)}$ has an irreducible summand of the form $S^{(u,v)}$ with $v\gs7$.
\end{propn}

\begin{pf}
The proof is very similar to the proof of Proposition \ref{ab0}. We need to show that there is a pair $u,v$ such that $S^{(u,v)}$ is irreducible, $v\gs7$, $v\equiv3\ppmod4$ and $\mbinom{u-v}{u-a}$ is odd.  The condition for $S^{(u,v)}$ to be irreducible is $u-v\equiv-1\pmod{2^{l(v)}}$.

Again, we need three cases.
\begin{description}
\item[\fbox{$a=8,10, n-9$ or $n-7$}]
In this case, take $v=7$ (so $u=n-7$).  Since $n\equiv5\ppmod8$, we get $u\equiv6\ppmod8$, which means that $u-v\equiv7\ppmod8$ (so $S^{(u,v)}$ is irreducible), and the binomial coefficients
\[
\binom{u-7}0,\binom{u-7}1,\binom{u-7}2, \binom{u-7}3
\]
are all odd, which means that $\mbinom{u-7}{u-a}$ will be odd.
\item[\fbox{$n\equiv5\ppmod{16}$, $12 \ls a \ls n-11$}]
In this case, let
\[
n'=\frac{n+5}2,\quad a'=
\begin{cases}
\mfrac{a+2}2&(a\equiv2\ppmod4)\\[5pt]
\mfrac{a+4}2&(a\equiv0\ppmod4).
\end{cases}
\]
Then $n',a'$ satisfy the conditions of the proposition, and $n'<n$. So by induction there is a pair $u',v'$ such that
\[
v'\equiv3\pmod4,\qquad v'\gs7,\qquad u-v\equiv-1\pmod{2^{l(v')}},\qquad \mbinom{u'-v'}{u'-a'}\equiv1\pmod2.
\]
Note that because $u'-v'$ is odd and $u'-a'$ is even, this also gives $\mbinom{u'-v'}{u'-a'+1}$ odd.

We let $u=2u'-2$ and $v=2v'-3$.  Then $u+v=n$, and we have
\[
u-v=2(u'-v')+1\equiv-1\ppmod{2^{l(v')+1}}
\]
and $l(v)\ls l(v')+1$.  So $S^{(u,v)}$ is irreducible.  Furthermore, $v\gs7$, $v\equiv3\ppmod4$ and
\[
\binom{u-v}{u-a}=\binom{2u'-2v'+1}{2u'-2a'(+2)}\equiv\binom{u'-v'}{u'-a'(+1)}\equiv1\pmod2,
\]
as required.
\item[\fbox{$n\equiv13\ppmod{16}$}]
In this case, let
\[
n'=\frac{n+13}2,\quad a'=
\begin{cases}
\mfrac{a+6}2&(a\equiv2\ppmod4)\\[5pt]
\mfrac{a+8}2&(a\equiv0\ppmod4).
\end{cases}
\]
Then $n',a'$ satisfy the conditions of the proposition, and $n'<n$. So by induction there is a pair $u',v'$ such that
\[
v'\equiv1\pmod4,\qquad v'\gs7,\qquad u-v\equiv-1\pmod{2^{l(v')}},\qquad\mbinom{u'-v'}{u'-a'}\equiv1\pmod2.
\]
Because $u'-v'$ is odd and $u'-a'$ is even, this also gives $\mbinom{u'-v'}{u'-a'+1}$ odd.

We let $u=2u'-6$ and $v=2v'-7$.  Then $u+v=n$, and we have
\[
u-v=2(u'-v')+1\equiv-1\pmod{2^{l(v')+1}},
\]
and $l(v)\ls l(v')+1$.  So $S^{(u,v)}$ is irreducible.  Furthermore, we have $v\equiv3\ppmod4$, $v\gs7$ and
\[
\binom{u-v}{u-a}=\binom{2u'-2v'+1}{2u'-2a'(+2)}\equiv\binom{u'-v'}{u'-a'(+1)}\equiv1\pmod2.\tag*{\raisebox{-10pt}{\qedhere}}
\]
\end{description}
\end{pf}

Now we can prove our main result.

\begin{pf}[Proof of Corollary \ref{maincor}]
Suppose we have a pair $a,b$ of positive even integers with $a\gs4$.  If $a\gs6$, $b\gs4$ and $a+b\equiv0\ppmod8$, then the result follows from Proposition \ref{ab0}. If $a\gs8$, $b\gs4$ and $a+b\equiv2\ppmod8$, then the result follows from Proposition \ref{ab1}.  If $a=6$ and $a+b\equiv2\ppmod8$, then by Theorem \ref{main} the Specht module $S^{(a+b,3)}$ is an irreducible summand of $S^\la$.

The second case of the corollary is precisely the condition for $S^{(a+b,3)}$ to be an irreducible summand of $S^\la$, while the third case is the condition for $S^{(a+b-4,5,2)}$ to be a summand.  So in any of the given cases, $S^\la$ certainly has an irreducible Specht module as a summand.  To complete the proof, it suffices to show that if $a=4$, $b=2$ or $a+b\equiv4$ or $6\ppmod8$, then the only possible Specht modules which can occur as irreducible summands of $S^\la$ are $S^{(a+b,3)}$ and $S^{(a+b-4,5,2)}$.

Suppose $S^\la$ has an irreducible summand $S^{(u,v)}$ with $v>3$. Then $v\equiv3\ppmod4$ and $u-v\equiv7\ppmod8$, which means that $u+v\equiv5\ppmod8$ and hence $a+b\equiv2\ppmod8$. Furthermore, $(u,v)\dom\la\reg$, which implies that $a\gs6$ and $b\gs4$.

Similarly, if $S^\la$ has an irreducible summand $S^{(u,v,2)}$ with $v>5$, then $v\equiv1\ppmod4$, $u-v\equiv7\ppmod8$ and hence $u+v+2\equiv3\ppmod8$, which gives $a+b\equiv0\ppmod8$. Furthermore, the fact that $(u,v,2)\dom\la\reg$ implies that $a\gs6$ and $b\gs4$.
\end{pf}

\section{Concluding remarks}\label{concsec}

The results in this paper do not give anything like a complete picture; this work is intended as a re-awakening of a long-dormant subject.  Given how small the first example of a decomposable Specht module in this paper is, it is surprising that it has taken thirty years for this example to be found.  We hope that this paper will be the start of a longer study of decomposable Specht modules.

We conclude the paper by making some speculations about decomposable Specht modules; these are based on calculations and observations, but we do not have enough evidence to make formal conjectures.

\subsection{Specht filtrations}

Our main results show that in certain cases summands of Specht modules are isomorphic to irreducible Specht modules.  In fact, reducible Specht modules can also occur as summands; for example, the first new decomposable Specht module $S^{(4,3,1^2)}$ found in this paper decomposes as $S^{(6,3)}\oplus S^{(4,3,2)}$, with the latter Specht module being reducible.

However, it is certainly not the case that every summand of a decomposable Specht module is isomorphic to a Specht module. But in the cases we have been able to calculate, every summand appears to have a filtration by Specht modules.  If this is true in general, it means that our main results are stronger, in that we have found all irreducible summands of Specht modules in our family.

In fact, we speculate that every Specht module has a filtration in which the factors are isomorphic to indecomposable Specht modules; this would imply in particular that every indecomposable summand has a Specht filtration.  This speculation is certainly true in the case of Specht modules labelled by hook partitions; this follows from \cite[\S2]{gm2}.

\subsection{$2$-quotient separated partitions}

In \cite[Definition 2.1]{jm}, James and Mathas make the following definition: a partition $\la$ is \emph{$2$-quotient separated} if it can be written in the form
\[
(c+2x_c,c-1+2x_{c-1},\dots,d+2x_d,d^{2y_d},(d-1)^{2y_{d-1}+1},\dots,1^{2y_1+1}),
\]
where $c+1\gs d\gs0$, $x_c\gs\cdots\gs x_d\gs0$ and $y_1,\dots,y_d\gs0$.  (Note that the definition includes the case $c=0$, where we have $\la=(2x_0)$ if $d=0$, or $(1^{2y_1})$ if $d=1$.)

Informally, the $2$-quotient separated condition means that the Young diagram of $\la$ can be decomposed as in the following diagram, where horizontal `dominoes' can appear in the first $c-d+1$ rows, and vertical `dominoes' can appear in the first $d$ columns.
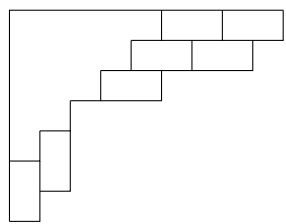
\[
\begin{tikzpicture}[scale=0.4]
\draw(0,0)--(0,7)--++(9,0);
\foreach \x in {0,1,2,3,4}\draw(\x,\x+2)--++(1,0)--++(0,1);
\foreach \x in {0,1}\draw(\x,\x)--++(1,0)--++(0,2);
\foreach \x in {0,1,2}\draw(\x+3,\x+4)--++(2,0)--++(0,1);
\foreach \x in {0,1}\draw(\x+6,\x+5)--++(2,0)--++(0,1);
\end{tikzpicture}
\]

The definition of a $2$-quotient separated partition was made as part of the study of decomposition numbers: for the Iwahori--Hecke algebra $\calh_{\bbc,-1}(\sss n)$, whose representation theory is very similar to that of $\sss n$ in characteristic $2$, the composition factors of a Specht module labelled by a $2$-quotient separated partition are known explicitly.  The reason we recall the definition here is that every known example of a decomposable Specht module is labelled by a $2$-quotient separated partition.  (Note that the partition $(a,3,1^b)$ considered in this paper is $2$-quotient separated precisely when $a$ and $b$ are even.)  It is interesting to speculate whether the $2$-quotient separated condition is necessary for a Specht module to be decomposable.

\end{document}